\documentclass{gtart}   
\usepackage{amssymb} 
\input gtmonout
\volumenumber{1}
\volumeyear{1998}
\volumename{The Epstein birthday schrift}
\pagenumbers{261}{293}
\received{2 June 1997}
\published{23 October 1998}
\papernumber{12}

\def\ra{\rightarrow}
\def\S{\Sigma}
\def\G{\Gamma}
\def\D{\Delta}

\def\<{\langle}
\def\>{\rangle}
\def\d{\partial}
\def\N{{\mathbb N}}
\def\Z{{\mathbb Z}}

\newtheorem{theorem}{Theorem}[section]

\newtheorem{lemma}[theorem]{Lemma}
\newtheorem{cor}[theorem]{Corollary}
\newtheorem{conj}[theorem]{Conjecture}
\theoremstyle{definition}
\newtheorem*{remark}{Remark}
\newtheorem*{remarks}{Remarks}

\begin{document}

\title{Minimal Seifert manifolds for higher ribbon knots}

\author{James Howie}

\address{Department of Mathematics, Heriot-Watt University\\ 
Edinburgh EH14 4AS, Scotland}
\email{J.Howie@hw.ac.uk}

\begin{abstract} 
We show that a group presented by a labelled oriented tree presentation
in which the tree has diameter at most three is an HNN extension of a finitely
presented group.   From results of Silver, it then follows that the 
corresponding higher dimensional ribbon knots admit minimal Seifert
manifolds.
\end{abstract}




\maketitlepage

\section{Introduction}\label{intro}

It is well known that every classical knot $k$ (knotted circle in
$S^3$) bounds a compact orientable surface, known as a {\it Seifert
surface} for the knot.   A Seifert surface $\S$ of minimal genus (among
all Seifert surfaces for the given knot $k$) is called {\it minimal},
and satisfies the following property: the inclusion-induced map
$\pi_1(\S\backslash k)\ra\pi_1(S^3\backslash k)$ is injective.

For a higher dimensional knot, or more generally a knotted (closed,
orientable) $n$--manifold $M$ in $S^{n+2}$, a {\it Seifert manifold} is
defined to be a compact, orientable $(n+1)$--manifold $W$ in $S^{n+2}$,
such that $\d W=M$.   A Seifert manifold $W$ for $M$ is defined to be
{\it minimal} if the inclusion-induced map $\pi_1(W\backslash
M)\ra\pi_1(S^{n+2}\backslash M)$ is injective.  In general, any $M$ will
always admit Seifert manifolds, but not necessarily minimal Seifert
manifolds.   For example, Silver \cite{S1} has shown that, for any 
$n\ge 3$, there exist $n$--knots in $S^{n+2}$ with no minimal
Seifert manifolds, and Maeda \cite{M} has constructed, for all $g\ge 1$,
a knotted surface of genus $g$ in
$S^4$ that has no minimal Seifert manifold.  Further examples of knotted tori
in $S^4$ without minimal Seifert manifolds are constructed by Silver \cite{S4}.

A theorem of Silver \cite{S2} says that, for $n\ge 3$, a knotted $n$--sphere $K$
in
$S^{n+2}$  has a minimal Seifert manifold if and only if its group
$G_K=\pi_1(S^{n+2}\backslash K)$ can be expressed as an HNN extension
with a {\it finitely presented} base group.   (It is standard that
any higher knot group can be expressed as an HNN extension with a
{\it finitely generated} base group.)

As Silver remarks, the proof of his theorem does not extend to the case
$n=2$.   However, it remains a {\it necessary} condition for the existence
of a minimal Seifert manifold that the group be an HNN extension with
finitely presented base group.   This applies also to knotted $n$--manifolds
in $S^{n+2}$, a fact which is used implicitly by Maeda in the result mentioned
above.   It remains an open
question whether every 2--knot in $S^4$ has a minimal Seifert manifold.
This seems unlikely, however.   For example Hillman \cite{Hi}, p. 139 shows that,
provided the 3--dimensional Poincar\'e Conjecture holds, there is an
infinite family of distinct 2--knots, all with the same group $G$,
such that the commutator subgroup of $G$ is finite of order 3; and at
most one of these knots can admit a minimal Seifert manifold.

In the present article we consider the case of higher dimensional
{\it ribbon knots}, for which the existence of minimal
Seifert manifolds is also an open question.   Indeed, as we shall
point out in the next section, higher ribbon knot groups are special
cases of {\it knot-like groups}, in the sense of Rapaport \cite{R},
and Silver \cite{S3} has conjectured that every finitely generated
HNN base for a knot-like group is finitely presented.   It would therefore
follow from Silver's conjecture (and his Theorem) that every higher
ribbon knot has a minimal Seifert manifold.

Now any higher ribbon knot group has a Wirtinger-like presentation
that can be encoded in the form of a {\it labelled oriented tree} (LOT) \cite{H1}.   Indeed the
LOT encodes not only a presentation for the knot group, but the 
complete homotopy type of the knot complement.   In \cite{H1} it was shown
that, if the diameter of the tree is at most 3, then the group is
locally indicable, and using this that the 2--complex model of the
associated Wirtinger presentation is aspherical.   A shorter proof of
this fact is given in \cite{I}, where it is shown that the presentation is in fact diagrammatically aspherical.

In the present
paper, we show that, under the same hypothesis on the diameter of the
tree, the group is an HNN extension with finitely presented base group,
and hence that the higher ribbon knot has a minimal Seifert manifold.

\begin{theorem}\label{main}\sl
Let $\G$ be a labelled oriented tree of diameter at most 3, and $G=G(\G)$
the corresponding group.   Then $G$ is an HNN extension with finitely
presented base group.
\end{theorem}

\begin{cor}\sl
Let $K$ be a ribbon $n$--knot in $S^{n+2}$, where $n\ge 3$, such that
the associated labelled oriented tree has diameter at most 3.   Then
$K$ admits a minimal Seifert manifold.
\end{cor}

The paper is arranged as follows.   In section \ref{defs} we recall
some basic definitions relating to LOTs and higher ribbon 
knots.   In section \ref{gps} we prove some preliminary results
about HNN bases for one-relator products of groups, which will
allow us to simplify the original problem.   In section \ref{reduction}
we reduce the problem to the study of {\it minimal} LOTs,
In section \ref{construction} we construct a finitely generated HNN base $B$ for
$G$, and describe a finite set of relators in these generators.  In
section \ref{struct} we prove some technical results about the structure of
these relations, which we apply in section \ref{kill} to complete the proof of 
 Theorem \ref{main} by
proving that this finite set is a set of defining relators for $B$.
We close, in section \ref{more}, with a geometric description of our
generators and relators for the HNN base, and a discussion of how this
might be used to generalise Theorem \ref{main}. 

\subsection*{Acknowledgements}

In the course of this work I have received useful comments and advice
from Nick Gilbert and from Dan Silver.   I am grateful to them for
their help.

\section{LOTs and higher ribbon knots}\label{defs}

A {\it labelled oriented tree} (LOT) is a tree $\G$, with vertex set $V=V(\G)$,
edge set $E=E(\G)$, and initial and terminal vertex maps
$\iota,\tau\co E\ra V$, together with an additional map $\lambda\co E\ra V$.
For any edge $e$ of $\G$, $\lambda(e)$ is called the {\it label}
of $e$.
In general, one can consider LOTs of any cardinality, but for the purposes
of the present paper, every LOT will be assumed to be finite.

\bigskip
To any LOT $\G$ we associate a presentation
$${\cal P}={\cal P}(\G):~~\<~V(\G)~|~\iota(e)\lambda(e)=\lambda(e)\tau(e)~\>$$
of a group $G=G(\G)$, and hence also a 2--complex $K=K(\G)$ modelled
on ${\cal P}$.   The 2--complex $K$ is a spine of a {\it ribbon disk complement}
$D^4\backslash k(D^2)$ \cite{H1}, that is the complement of an embedded 2--disk in
$D^4$, such that the radial function on $D^4$ composed with the embedding
$k$ is a Morse function on $D^2$ with no local maximum.
Conversely, any ribbon disk complement has a 2--dimensional spine of the
form $K(\G)$ for some LOT $\G$.

By doubling a ribbon disk, we obtain a ribbon 2--knot in $S^4$, and by
successively spinning we can obtain ribbon $n$--knots in $S^{n+2}$ for all
$n\ge 2$.   In each case the group of the knot is isomorphic to
the fundamental group of the ribbon disk complement that we started with.
Conversely, every ribbon $n$--knot (for $n\ge 2$) can be constructed this way,
so that higher ribbon knot groups and LOT groups are precisely the same
thing.

\bigskip
Recall \cite{R} that a group $G$ is {\it knot-like} if it has a finite
presentation with deficiency 1 (in other words, one more generator than
defining relator), and infinite cyclic abelianisation.   It is clear that
every LOT group has these properties, so LOT groups are special cases of knot-like groups.

\bigskip
The {\it diameter} of a finite connected graph $\G$ is the maximum distance
between two vertices of $\G$, in the edge-path-length metric.  A key factor
in our situation is the special nature of trees of diameter 3 or less.
For any LOT $\G$ of diameter 0 or 1, it is easy to see that $G(\G)$ is infinite
cyclic, so such LOTs are of little interest.

\begin{remark}
Every tree of diameter 2 has a single non-extremal vertex.   Every tree
of diameter 3 has precisely 2 non-extremal vertices.
\end{remark}

\bigskip
We recall from \cite{H1} that a LOT $\G$ is {\it reduced} if:

\begin{enumerate}
\item for all $e\in E$, $\iota(e)\ne\lambda(e)\ne\tau(e)$;
\item for all $e_1\ne e_2\in E$, if $\lambda(e_1)=\lambda(e_2)$ then
$\iota(e_1)\ne\iota(e_2)$ and $\tau(e_1)\ne\tau(e_2)$;
\item every vertex of degree 1 in
$\G$ occurs as a label of some edge of $\G$.
\end{enumerate}

For every LOT $\G$ there is a reduced LOT $\G'$ with the same group as $\G$,
and the same or smaller diameter, 
so we may also restrict our attention to reduced LOTs.

A subgraph $\G'$ of a LOT $\G$ is {\it admissible} if $\lambda(e)\in V(\G')$
for all $e\in E(\G')$.   If $\G'$ is connected and admissible, then it
is also a LOT.
A LOT is {\it minimal} if every connected admissible subgraph consists
only of a single vertex.

\bigskip
If $\G$ is a LOT and $A\subseteq V(\G)$, we define the {\it span} of $A$
(in $\G$) to be the smallest subgraph $\G'$ of $G$ such that:
\begin{enumerate}
\item $A\subseteq V(\G')$; and
\item if $e\in E(\G)$ with $\lambda(e)\in V(\G')$ and at least one 
of $\iota(e)$, $\tau(e)$ belongs to $V(\G')$, then $e\in E(\G')$.
\end{enumerate}

We write ${\rm span}(A)$ for the span of $A$, and say that $A$ {\it spans},
or {\it generates} $\G'$ if $\G'={\rm span}(A)$.   The following is essentially
Proposition 4.2 of \cite{H1}.

\begin{lemma}\sl
If $\G$ is a LOT spanned by $A$, then ${\cal P}(\G)$ is Andrews--Curtis
equivalent to a presentation with generating set $A$.   If $\G'$ is an
admissible subgraph of $\G$ with $V(\G')\subseteq A$, then the presentation
may be chosen to contain ${\cal P}(\G')$, and the Andrews--Curtis moves can
be taken relative to ${\cal P}(\G')$.
\end{lemma}

\begin{cor}\label{2-sp}\sl
If $\G$ is a LOT spanned by two vertices, then $G(\G)$ is a torsion-free
one-relator group.
\end{cor}

\begin{proof}
Let $A$ be a set of two vertices spanning $\G$.   Then ${\cal P}(\G)$ is
Andrews--Curtis equivalent to a presentation $\<A|R\>$.   Since ${\cal P}(\G)$
has deficiency 1, the same is true of the
equivalent presentation $\<A|R\>$.   In other words $|R|=1$, and $G(\G)$ is
a one-relator group.   But the abelianisation $G^{ab}$ of $G$ is infinite cyclic, so the relator $r\in R$ cannot be a proper power, and so $G$ is torsion-free.
\end{proof}

We will require the following generalisation of Corollary \ref{2-sp}.
Recall that a {\it one-relator product} of two groups $A,B$ is the
quotient of the free product $A*B$ by the normal closure of a single
word $w$, called the {\it relator}.

\begin{cor}\label{+1-sp}\sl
If $\G$ is a LOT spanned by $V(\G')\cup\{x\}$, where $\G'$ is an admissible
subgraph of $\G$ and $x$ is a vertex in $V(\G)\backslash V(\G')$, then
$G(\G)$ is a one-relator product of $G(\G')$ and $\Z$, where the relator is
not a proper power.
\end{cor}

\begin{proof}
Let $A=V(\G')\cup\{x\}$ and apply the Theorem.   Then ${\cal P}(\G)$ is
equivalent, relative to ${\cal P}(\G')$, to a presentation ${\cal Q}$
with generating set $A$ and containing ${\cal P}(\G')$.   Now each of
${\cal P}(\G)$, ${\cal P}(\G')$ and ${\cal Q}$ has deficiency 1. 
  Moreover, ${\cal Q}$ has one
more generator than ${\cal P}(\G')$, so ${\cal Q}$ also has one
more defining relator than ${\cal P}(\G')$.   It follows that $G(\G)$
is a one relator product of $G(\G')$ with the infinite cyclic group
$\<x\>$.   Finally, since the abelianisations of $G(\G)$, $G(\G')$ and
$\<x\>$ are all infinite cyclic, it follows that the relator cannot be
a proper power.
\end{proof}

\section{One-relator groups and one-relator products}\label{gps}

The following result is merely a summary of some well-known properties of one-relator
groups, which have useful applications to our situation.
Recall that a group $G$ is {\it locally indicable} if, for every
nontrivial, finitely generated subgroup $H$ of $G$, there exists an
epimorphism $H\ra\Z$.

\begin{theorem}\label{1-rel-gp}\sl
Let $G$ be a finitely generated one-relator group.   Then
\begin{enumerate}
\item $G$ is either a finite cyclic group, or an HNN extension of a 
finitely presented, one-relator group (with shorter defining relator);
\item if the defining relator of $G$ is not a proper power, then $G$
is locally indicable.
\end{enumerate}
\end{theorem}

\begin{proof}
See \cite{Mo} and \cite{Br} respectively.
\end{proof}

In order to complete the process of reducing ourselves to a simple
special case, we require a generalisation of the above theorem to
one-relator products.   Suppose that $A$ and $B$ are locally indicable
groups, and $N=N(w)$ is the normal closure in $A*B$ of a cyclically
reduced word $w$ of length at least 2 that is not a proper power.   
Then the one-relator product $G=(A*B)/N$ is known \cite{H2}
to be locally indicable.   We show also that $G$ has a finitely presented
HNN base, provided that $A$ and $B$ also have this property.

\begin{theorem}\sl
Let $G=(A*B)/N(w)$ be a one-relator product of two finitely presented,
locally indicable groups $A$ and $B$, each of which has a finitely
presented HNN base.   Suppose also that $G^{ab}$ is infinite cyclic,
with each of the natural maps $A^{ab}\ra G^{ab}$ and $B^{ab}\ra G^{ab}$
an isomorphism.   Then $G$ is a finitely presented,
locally indicable group with a finitely presented HNN base.
\end{theorem}

\begin{remark}
The condition on $G^{ab}$ in this theorem is unnecessary for the proof
that $G$ has a finitely presented HNN base.   It can
be removed at the expense of a less straightforward proof.   However
the condition does hold for all the groups that we are considering
in this paper, so there is no loss of generality for us in imposing
that condition.   The condition also ensures that $w$ cannot be a proper
power, so that $G$ is locally indicable by the results of \cite{H2}.
\end{remark}

\begin{proof}
A presentation for $G$ can be obtained by taking the disjoint union
of finite presentations for $A$ and for $B$, and imposing the
single additional relation $w=1$.   Hence $G$ is finitely presented.
As pointed out in the remark above, $w$ cannot be a proper power,
so $G$ is locally indicable by \cite{H2}.   It remains only to
prove that $G$ has a finitely presented HNN base.

Let
$$A=\<A_0,a|a^{-1}ga=\alpha(g)~(g\in A_1)\>$$
and
$$B=\<B_0,b|b^{-1}hb=\beta(h)~(h\in B_1)\>$$
be HNN presentations for $A$ and $B$ with finitely presented bases $A_0$
and $B_0$ respectively.   Since $A$ and $B$ are finitely presented,
it follows also that the associated subgroups $A_1$ and $B_1$ are finitely
generated.

The commutator subgroup $G'$ of $G$ can be expressed in the form
$$(A'*B'*\<~c_n~(n\in\N)\>)/N(\{w_n~(n\in\N)\}),$$
where $c_n=a^{n+1}b^{-1}a^{-n}$ and $w_n=a^{-n}wa^n$.

Now $A'$ is an infinite stem product

\begin{center}
\begin{tabular}{ccccccc}
$\cdots$ &$(a^{-1}A_0a)$ &$*$ &$A_0$ &$*$ &$(aA_0a^{-1})$ &$\cdots$ \\
 & &${(a^{-1}A_1a)}$ & &${A_1}$ & &
\end{tabular}
\end{center}
Since $A_0$ is finitely presented and $A_1$ is finitely generated,
the subgroup
\begin{center}
\begin{tabular}{ccccc}
$(a^{-k}A_0a^k)$ &$*$ &$\cdots\cdots$ &$*$ &$(a^kA_0a^{-k})$ \\
 &$(a^{-k}A_1a^k)$ & &$(a^{k-1}A_1a^{1-k})$ &
\end{tabular}
\end{center}
is finitely presented for each $k$.   Moreover it is also an HNN base
for $A$.   Replacing $A_0$ by this subgroup, for any sufficiently large $k$, we may assume that $w_0\in A_0*B'*\<~c_n~(n\in\N)\>$.

Similarly, possibly after replacing $B_0$ by a sufficiently large finitely presented
HNN base for $B$, we may assume that $w_0\in A_0*B_0*\<~c_n~(n\in\N)\>$.   Now
let $\mu$ and $\nu$ be the least and greatest indices $i$ such that
$c_i$ occurs in $w_0$.   (Note that at least one $c_i$ occurs in $w_0$,
for otherwise $w_0\in A_0 * B_0$, so $w\in A'*B'$, whence
$G^{ab}\cong A^{ab}\times B^{ab}\not\cong\Z$, a contradiction.)
Define $G_0=(A_0*B_0*\<c_\mu,\dots,c_\nu\>)/N(w_0)$ and 
$G_1=A_0*B_0*\<c_\mu,\dots,c_{\nu-1}\>$, and observe that $G_0$ is
a finitely presented HNN base for $G$, with associated subgroup
$G_1$.
\end{proof}

\section{Reduction of the problem}\label{reduction}

Recall from section \ref{defs} that a LOT $\G$ is {\it minimal} if it contains no admissible subtree with more
than one vertex.   In this section we reduce the proof of the main
theorem to the case of a minimal LOT of diameter 3, using the results
of section \ref{gps}.   The key point is that a non-minimal LOT can
be obtained from a minimal admissible subtree by successively expanding
to the span of the existing tree with one extra vertex.   By Corollary
\ref{+1-sp}, this
construction corresponds at the group level to taking a one-relator 
product of a given group with an infinite cyclic group.

\begin{lemma}\label{reduce}\sl
Let $\G$ be a LOT of diameter at most 3, containing a proper admissible
subtree with more than one vertex.   Then there is such an admissible
subtree $\G'$ and a vertex $x\in V(\G)\backslash V(\G')$ such that $\G$
is spanned by $V(\G')\cup\{x\}$.
\end{lemma}

\begin{proof}
Suppose first that some extremal vertex $x$ of $\G$ does not occur as
a label of any edge of $\G$.   In this case we take $\G'$ to consist of $\G$ with the vertex $x$ and the edge incident to $x$ removed.   Clearly $\G'$ is
connected, so a subtree of $\G$.  Since $x$
is not the label of any edge in $E(\G')$, it follows that $\G'$ is 
admissible.   Moreover $\G$ is spanned by $V(\G)=V(\G')\cup\{x\}$, as
required.

We may therefore assume that every extremal vertex of $\G$ occurs at
least once as the label of an edge of $\G$.

Next suppose that $\G$ has a proper admissible subtree that contains all
the non-extremal vertices of $\G$.   Let $\G'$ be a maximal such
admissible subtree.   The vertices in $V(\G)\backslash V(\G')$
are all extremal in $\G$, so occur as labels of edges of $\G$.
But since $\G'$ is admissible, no such vertex can be a label of an
edge of $\G'$.   Since the finite sets $V(\G)\backslash V(\G')$
and $E(\G)\backslash E(\G')$ have the same cardinality, it follows that
each vertex in $V(\G)\backslash V(\G')$ is the label of precisely
one edge in $E(\G)\backslash E(\G')$.   In turn, this edge has precisely
one endpoint in $V(\G)\backslash V(\G')$, so we can define a permutation
$\sigma$ on $V(\G)\backslash V(\G')$ by defining $\sigma(x)$ to
be the extremal endpoint of the unique edge labelled $x$, for all 
$x\in V(\G)\backslash V(\G')$.   Now fix some vertex 
$x\in V(\G)\backslash V(\G')$, let $t$
be the size of the orbit of $\sigma$ that contains $x$, and
define $x_i=\sigma^i(x)$, $i=1,\dots,t$.   Now 
$\Delta={\rm span}(V(\G')\cup\{x\})$
contains the vertex $x=x_t$, together with any non-extremal vertex
of $\G$.  Hence $\Delta$ contains the edge labelled $x_t$, and hence
its endpoint $x_1$.   Similarly $\Delta$ contains $x_2,\dots,x_{t-1}$,
as well as the edges labelled $x_1,\dots, x_{t-1}$.   On the other hand,
The vertices $x_1,\dots, x_t$, the edges labelled by them, and
the vertices and edges of $\G'$ together form an admissible
subtree of $\G$, which by maximality of $\G'$ must be the whole of $\G$.
Hence $\Delta=\G$, in other words $\G$ is spanned by $V(\G')\cup\{x\}$.

Finally, suppose that no proper 
admissible subtree of $\G$ contains all the
non-extremal vertices of $\G$.  In particular, $\G$ must have more than
one non-extremal vertex, so has diameter 3.   By hypothesis, there is
a proper admissible subtree $\G'$ of $\G$ that contains more than one
vertex.   Hence $\G'$ contains precisely one of the two nonextremal
vertices of $\G$, say $u$.   As an abstract graph, $\G$ is the union
of $\G'$ with another tree $\G''$, such that $\G'\cap\G''=\{u\}$.   Note
that $\G''$ contains both of the non-extremal vertices of $\G$, so
cannot be an admissible subtree, by hypothesis.   Hence at least one
edge $f$ of $\G''$ is labelled by a vertex $a$ of $\G'$ (other than $u$).
Let $e$ be the edge of $\G$ that joins the two non-extremal vertices 
$u,v$, and let $\Delta={\rm span}(V(\G')\cup\{\lambda(e)\})$.   Then
$\Delta$ contains $\G'$ and the edge $e$, and hence $v$, and hence the edge $f$.   Each extremal vertex of $\Delta$ is the label of an edge
of $\G$, and hence of $\Delta$, since $\Delta$ contains at least one
endpoint (namely $u$ or $v$) of every edge of $\G$.   Moreover there are
$|E(\G')|+1$ edges of $\Delta$ labelled by the $|V(\G')|=|E(\G')|+1$
vertices of $\G'$, so an easy counting argument shows that there must be
at least $|V(\Delta)|-1$ edges in $\Delta$.   In other words $\Delta$
is a tree, so the whole of $\G$.   In other words $\G$ is spanned by
$V(\G')\cup\{\lambda(e)\}$.
\end{proof}

\begin{remark}
If $\G$ is a minimal LOT of diameter 2, 
then the above argument still applies (to the subtree
consisting of only the unique non-extremal vertex).   In this
case we see that the permutation $\sigma$ is transitive, and that
$\G$ is spanned by two vertices.
\end{remark}

\begin{lemma}\label{span}\sl
Let $\G$ be a minimal LOT of diameter 3, and let $u,v$ be the
two non-extremal vertices of $\G$.  Then one of the following holds:
\begin{enumerate}
\item One of $u,v$ is a label in $\G$, and $\G$ is spanned by $\{u,v\}$;
\item Some vertex $a$ occurs twice as a label in $\G$, and $\G$ is
spanned by $\{a,u,v\}$.
\end{enumerate}
\end{lemma}

\begin{proof}
By minimality of $\G$, every extremal vertex of $\G$ occurs as a label.
There are $|V|-2$ extremal vertices, and $|V|-1$ edges, so either
one of $u,v$ occurs as a label or some unique extremal vertex $a$
occurs twice as a label.   Note that every edge of $\G$ is incident
to at least one of $u,v$, so if $u,v\in A\subset V$ then every edge
labelled by a vertex of ${\rm span }(A)$ is an edge of ${\rm span }(A)$.
\begin{enumerate}
\item Suppose that $u$ occurs as a label, and let $\G'={\rm span }(\{u,v\})$.
If $\G'$ has $k+2$ vertices $u,v,x_1,\dots,x_k$, then $x_1,\dots,x_k$
are all extremal in $\G$, so each of $u,x_1,\dots,x_k$ is a label of
an edge of $\G$, which must therefore be an edge of $\G'$.   Hence
$\G'$ has at least $k-1$ edges, so is connected.   By minimality of
$\G$ we have $\G=\G'={\rm span }(\{u,v\})$.
\item Suppose that an extremal vertex $a$ appears twice as a label,
and let $\G'={\rm span }(\{a,u,v\})$.   If $\G'$ has $k+3$ vertices
$a,u,v,x_1,\dots,x_k$, then each of $x_1,\dots,x_k$ is extremal, so
the label of an edge of $\G$, while $a$ is the label of 2 edges of $\G$.
Each of these $k+2$ edges is an edge of $\G'$, so $\G'$ is connected,
and by minimality again we have  $\G=\G'={\rm span }(\{a,u,v\})$.
\end{enumerate}

\vspace{-1.5\baselineskip}\end{proof}

\begin{cor}\sl
If $\G$ is either a minimal LOT of diameter 2, or a minimal LOT of
diameter 3 in which no vertex occurs twice as a label, then $G(\G)$
is a locally indicable group with a finitely presented HNN base.
\end{cor}

\begin{proof}
By Lemma \ref{span} or the remark following Lemma \ref{reduce}, $\G$ is spanned
by two vertices.   Hence $G=G(\G)$ is a 2--generator, one-relator
group.  Since $G^{ab}$ is infinite cyclic, $G$ is not finite, and the relator of $G$ cannot be a proper power.   The result follows
immediately from Theorem \ref{1-rel-gp}.
\end{proof}

Using the above results, we can reduce our problem to the case of
a minimal LOT of diameter 3 that is not spanned by two vertices.
In particular, some extremal vertex must occur twice as a label.

\begin{cor}\sl
If the group of every reduced, minimal LOT of diameter 3 which is
not spanned by two vertices is locally indicable
with finitely presented HNN base, then the same is true for every
LOT of diameter 3 or less.
\end{cor}

Recall \cite{H1} that the {\it initial graph} $I(\G)$ of $\G$ is
the graph with the same vertex and edge sets as $\G$, but with
incidence maps $\iota,\lambda$.   Similarly the {\it terminal graph}
$T(\G)$ of $\G$ has the same vertex and edges sets as $\G$, but incidence
maps $\lambda,\tau$.   It was shown in \cite{H1} that the commutator
subgroup of $G(\G)$ is locally free if either $I(\G)$ or $T(\G)$ is
connected.   (If $I(\G)$ and $T(\G)$ are both connected, then $G(\G)'$
is free of finite rank.)   In particular, any finitely generated
HNN base for $G(\G)$ is free, so automatically finitely presented.

Hence we can concentrate attention on the case of a minimal LOT
$\G$ of diameter 3, not spanned by any two of its vertices, such
that neither $I(\G)$ nor $T(\G)$ is connected.   Our next result
gives a detailed description of the structure of $I(\G)$.   In
particular it will show us that $I(\G)$ has precisely two connected
components, one containing each of the nonextremal vertices of
$\G$.   A similar statement holds for $T(\G)$.

\begin{lemma}\label{I}\sl
Let $\G$ be a minimal LOT of diameter 3, with nonextremal vertices $u$
and $v$, and an extremal vertex $a$ that occurs twice as a label of
edges of $\G$.   Then:
\begin{enumerate}
\item $u$ and $v$ are sources in $I(\G)$;
\item no vertex other than $u$ or $v$ is the initial vertex of more than
one edge of $I(\G)$;
\item $a$ is the terminal vertex of precisely two edges of $I(\G)$;
\item each vertex other than $a,u,v$ is the terminal vertex of precisely
one edge of $I(\G)$;
\item any directed cycle in $I(\G)$ contains $a$;
\item each component of $I(\G)$ contains at least one of $u,v$;
\item $I(\G)$ has at most two connected components.
\end{enumerate}
\end{lemma}

\begin{proof}

\vspace{-.5\baselineskip}\begin{enumerate}
\item Since $\lambda(e)\ne u$ for all $e\in E(\G)$, $u$ is not the terminal
vertex of any edge in $I(\G)$, in other words $u$ is a source. 
Similarly $v$ is a source in $I(\G)$.
\item
Any vertex $x$ of $\G$, with the exception of $u$ and $v$, is extremal
in $\G$, so the initial vertex of at most one edge of $\G$.   Hence
$x$ is also the initial vertex of at most one edge in $I(\G)$.
\item
$a=\lambda(e)$ for precisely two edges $e\in E(\G)$.
\item
If $x\in V(\G)\backslash\{a,u,v\}$ then $x=\lambda(e)$ for precisely
one edge $e\in E(\G)$. 
\item
Suppose $(e_1,e_2,\dots,e_n)$ is a directed cycle in $I(\G)$.
Then there are vertices $x_1,\dots,x_n\in V(\G)$ with $x_i=\iota(e_i)$
for all $i$, $\lambda(e_i)=x_{i+1}$ for $i<n$, and $\lambda(e_n)=x_1$.
Now each $x_i$ is extremal since it occurs as a label.   If no
$x_i$ is equal to $a$ then we can remove the vertices $x_1,\dots,x_n$
and the edges $e_1,e_2,\dots,e_n$ from $\G$ to form a connected,
admissible subgraph $\G'$ that contains at least three vertices ($a,u,v$).
This contradicts the minimality of $\G$, and so $x_i=a$ for some $i$,
as claimed.
\item
By (iv) if $x\not\in\{a,u,v\}$ then $x$ is the terminal vertex in
$I(\G)$ of a unique edge.   If the initial vertex of this edge is not one
of $a,u,v$ then it also is the terminal vertex of a unique edge.
Continuing in this way, we can construct a directed path that ends at
$x$, and either begins at
one of $a,u,v$ or contains a cycle.   By (v) any directed
cycle contains $a$, so in any case we have a directed path from
one of $a,u,v$ to $x$.   It suffices therefore to find a path in
$I(\G)$ from $u$ or $v$ to $a$.   But $a$ is the terminal vertex in 
$I(\G)$ of precisely two edges, with initial vertices $x_1$ and $x_2$
say.   Now apply the above argument to each of $x_1,x_2$.   If there
is a path from $u$ or $v$ to $x_1$ or $x_2$ then we are done.
Otherwise there are directed paths from $a$ to each of $x_1,x_2$.
Neither $u$ nor $v$ can belong to these paths, since they are sources
in $I(\G)$.   But then from (ii) it follows that there is at most one 
directed path of any given length beginning at $a$, whence $x_1=x_2$,
a contradiction.   Hence there is a directed path in $I(\G)$ from
$u$ or $v$ to $a$, as claimed.
\item
This follows immediately from (vi).
\end{enumerate}

\vspace{-1.5\baselineskip}\end{proof}

A similar result holds for $T(\G)$.

\begin{lemma}\label{T}\sl
Let $\G$ be a minimal LOT of diameter 3, with nonextremal vertices $u$
and $v$, and an extremal vertex $a$ that occurs twice as a label of
edges of $\G$.   Then:
\begin{enumerate}
\item $u$ and $v$ are sinks in $T(\G)$;
\item no vertex other than $u$ or $v$ is the terminal vertex of more than
one edge of $T(\G)$;
\item $a$ is the initial vertex of precisely two edges of $T(\G)$;
\item each vertex other than $a,u,v$ is the initial vertex of precisely
one edge of $T(\G)$;
\item any directed cycle in $T(\G)$ contains $a$;
\item each component of $T(\G)$ contains at least one of $u,v$;
\item $T(\G)$ has at most two connected components.
\end{enumerate}
\end{lemma}

\begin{cor}\label{IT}\sl
Suppose that $\G$ is a reduced, minimal LOT of diameter 3,
which is not spanned by two vertices, and such that neither $I(\G)$
nor $T(\G)$ is connected.   Then
\begin{enumerate}
\item
There is a unique extremal vertex $a$ of $\G$ that is the label of
two distinct edges of $\G$.   One of these edges has an extremal
initial vertex, and the other has an extremal terminal vertex.
\item $I(\G)$ has precisely two connected components, each containing
one of the two nonextremal vertices $u,v$ of $\G$. 
\item  There is a unique cycle
in $I(\G)$, which is either a directed cycle containing $a$, or consists
of two directed paths (one of length 1, the other of length at least 2),
from $u$ or $v$ to $a$.  
\item $T(\G)$ has precisely two connected components, each containing
one of the two nonextremal vertices $u,v$ of $\G$. 
\item  There is a unique cycle
in $T(\G)$, which is either a directed cycle containing $a$, or consists
of two directed paths (one of length 1, the other of length at least 2),
from $a$ to $u$ or $v$.  
\item The cycles in $I(\G)$ and $T(\G)$ are not both directed.
\end{enumerate}
\end{cor}

\begin{proof}

\vspace{-.5\baselineskip}\begin{enumerate}
\item
We already know that there is an extremal vertex $a$ occurring twice
as a label, by Lemma \ref{span}, since otherwise $\G$ can be spanned by two vertices.
We also know that $a$ is unique, since every extremal vertex occurs at
least once as a label.   Now suppose that neither of the edges labelled
$a$ has extremal initial vertex.   The initial vertices of these two
edges must be distinct, since $\G$ is reduced, and so must
be the two nonextremal vertices $u,v$ of $\G$.   But then there
are edges of $I(\G)$ from both $u$ and $v$ to $a$. Hence $u$ and $v$
belong to the same connected component of $I(\G)$.   By Lemma \ref{I}, (vi)
it follows that $I(\G)$ is connected, a contradiction. 
A similar contradiction arises if neither edge has an extremal
terminal vertex.
\item
This is just a restatement of Lemma \ref{I}, (vi), together with the
hypothesis that $I(\G)$ is not connected.
\item
Since $I(\G)$ has the same vertex and edge sets as $\G$, it has the same
euler characteristic, namely 1.   Since $I(\G)$ has two components,
it follows that $H_1(\G)\cong\Z$, so there is a unique cycle in $I(\G)$.
If this cycle is directed, then it must contain $a$, by Lemma \ref{I}, (v).
Otherwise it must contain at least two vertices at which the orientation
of the edges of the cycle changes.   This is possible only at a vertex
which is either the initial vertex of at least two edges or the
terminal vertex of at least two edges, and by Lemma \ref{I} the only such
vertices are $a,u,v$.   Let us assume that $a$ is in the same component
of $I(\G)$ as $u$.   Then the cycle must contain both $a$ and $u$, and 
indeed must consist of two directed paths from $u$ to $a$.   
By uniqueness of the cycle (or directly from Lemma \ref{I}), we see that
there only two directed paths in $I(\G)$ from $u$ to $a$.   Moreover,
precisely one of these paths is of length 1, since precisely one
of the edges of $\G$ labelled $a$ has a nonextremal initial vertex.
\item Similar to (ii).
\item Similar to (iii).
\item
If the cycle in $I(\G)$ is directed, then there is an edge of $I(\G)$
with initial vertex $a$, and so also there is an edge of $\G$ with initial
vertex $a$.   Similarly, if the cycle in $T(\G)$ is directed, then there
is an edge of $\G$ with terminal vertex $a$.   Since $a$ is extremal
in $\G$, these cannot both occur.
\end{enumerate}

\vspace{-1.4\baselineskip}\end{proof}

\section{Construction of the HNN base}\label{construction}

In this section, we construct a presentation of a group that
will turn out to be an HNN base for $G$.   As a first step,
we fix names for the various vertices of $\G$.   Throughout
we make the following assumptions:

\begin{itemize}
\item $\G$ is a minimal LOT of diameter 3, which cannot be spanned by
fewer than three vertices.
\item The non-extremal vertices of $\G$ are $u$ and $v$.
\item The unique vertex of $\G$ that appears twice as a label is $a$.
\item Of the edges labelled $a$, one has its initial vertex in $\{u,v\}$
and its terminal vertex extremal, while the other has its initial vertex
extremal and its terminal vertex in $\{u,v\}$.
\item Neither $I(\G)$ nor $T(\G)$ is connected.
\end{itemize}

We know from Lemma \ref{span} that $\G$ is then spanned by $\{a,u,v\}$.  Let
$\D$ denote the subtree of $\G$ whose vertex set is $\{a,u,v\}$.
We give inductive definitions of 
two sequences $\{b_1,b_2,\dots,b_P\}$ and $\{c_1,c_2,\dots,c_Q\}$
of vertices of $\G$, and two sequences $\{e_0,\dots,e_P\}$,
$\{f_0,\dots,f_Q\}$ of edges of $\G$ as follows.   

Define $e_0$ to be the edge of $\G$ whose label is $a$ and whose terminal vertex
is in $\{u,v\}$.   For $i\ge 0$, assume inductively that
$e_i$ has been defined. If $e_i$ is an edge of $\D$, then we define $P=i$ and stop
the construction of the sequences $\{b_1,b_2,\dots,b_P\}$ and
$\{e_0,\dots,e_P\}$.
Otherwise $e_i$ joins one of $\{u,v\}$ to an extremal vertex other than
$a$, and we define $b_{i+1}$ to be that extremal vertex, and 
$e_{i+1}$ to be the unique edge of $\G$ labelled 
$b_{i+1}$.  

Similarly, define $f_0$ to be the edge of $\G$ whose label is $a$ and whose
initial vertex is in $\{u,v\}$.
For $i\ge 0$, assume inductively that $f_i$ has been defined.   If $f_i$ is an edge of $\D$, then we 
define $Q=i$ and stop the construction of the sequences 
$\{c_1,c_2,\dots,c_Q\}$ and $\{f_0,\dots,f_Q\}$.
Otherwise $f_i$ joins one of $\{u,v\}$ to an extremal vertex other than
$a$, and we define $c_{i+1}$ to be that extremal vertex, and $f_{i+1}$
to be the unique edge labelled by $c_{i+1}$.

Note that the $P+Q+3$ vertices $\{u,v,a,b_1,\dots,b_P,c_1,\dots,c_Q\}$
and the $P+Q+2$ edges $\{e_0,\dots,e_P,f_0,\dots,f_Q\}$ together form
an admissible subgraph of $\G$, which has euler characteristic 1 and hence
is connected, and hence by minimality of $\G$ must be the whole of $\G$.
In other words 
$$V=V(\G)=\{u,v,a,b_1,\dots,b_P,c_1,\dots,c_Q\},$$ and
$$E=E(\G)=\{e_0,\dots,e_P,f_0,\dots,f_Q\}.$$

We also introduce the following notation.   For $i=1,\dots,P$, $x_i$
denotes the unique non-extremal vertex of $\G$ (ie $x_i\in\{u,v\}$)
incident with the edge $e_{i-1}$.   For $i=1,\dots,Q$, $y_i$ 
denotes the unique non-extremal vertex of $\G$ 
incident with the edge $f_{i-1}$.   In other words, $x_i$ is the vertex
adjacent to $b_i$ in $\G$, and $y_i$ is the vertex adjacent to $c_i$.

\begin{lemma}\label{key}\sl

\begin{enumerate}
\item If $x_2=\dots=x_P=u$, then $x_1=v$ and $e_P$ is incident at $v$.
\item If $x_2=\dots=x_P=v$, then $x_1=u$ and $e_P$ is incident at $u$.
\item If $y_2=\dots=x_Q=u$, then $y_1=v$ and $f_Q$ is incident at $v$.
\item If $y_2=\dots=y_Q=v$, then $y_1=u$ and $f_Q$ is incident at $u$.
\end{enumerate}
\end{lemma}

\begin{proof}
We prove (i).   The other proofs are similar.

Suppose first that $x_1=x_2=\dots=x_P=u$, and consider the subgraph
$\G_0={\rm span}\{a,u\}$ of $\G$.   Since $\lambda(e_0)=a$ and $e_0$
is incident to $u$, we have $e_0\in E(\G_0)$, and since $b_1$ is an
endpoint of $e_0$ we have $b_1\in V(\G_0)$.   Similarly $e_1\in E(\G_0)$
and $b_2\in V(\G_0)$, and so on, until $e_P\in E(\G_0)$.
If $e_P$ is incident with $v$, then $v\in V(\G_0)$, and since $\G$ is spanned by $\{a,u,v\}$ it follows that $\G=\G_0$ is spanned by $\{a,u\}$,
a contradiction.   Otherwise, $e_P$ joins $a$ to $u$, in which case
the vertices $a,u,p_1,\dots,b_P$ and the edges $e_0,\dots,e_P$ form
an admissible subtree of $\G$ of diameter two, which again is a
contradiction.

Now suppose that $x_1=v$ and $x_2=\dots=x_P=u$, and let $\G_0={\rm span}\{b_1,u\}$.
Arguing as above, we see that $\G_0$ contains the edges $e_1,\dots,e_{P-1}$
and the vertices $u,b_1,\dots,b_P$.   If $e_P$ is not incident at $v$,
then it joins $u$ to $a$, so $e_P$ and $a$ also belong to $\G_0$.
But then $e_0$ joins $b_1$ to $v$ and has
label $a$, so we also have $v\in V(\G_0)$.  Hence $\G=\G_0$ since $\G$ is spanned by $\{a,u,v\}$,
and so $\G$ is spanned by $\{b_1,u\}$, a contradiction.
\end{proof}

\bigskip
We next subdivide each of the sequences $\{b_i\}$, $\{c_i\}$ into
two subsequences, depending on the orientation of the edges labelled
by these vertices.   Specifically, let:

\begin{itemize}
\item $p(1),\dots,p(s)$ be the sequence,
in ascending order, of integers $i$ such that $0<i\le P$ and $b_i=\tau(e_{i-1})$;
\item $p'(1),\dots,p'(s')$ be the sequence,
in ascending order, of integers $i$ such that $0<i\le P$ and $b_i=\iota(e_{i-1})$;
\item $q(1),\dots,q(t)$ be the sequence,
in ascending order, of integers $i$ such that $0<i\le Q$ and $c_i=\iota(f_{i-1})$; and
\item $q'(1),\dots,q'(t')$ be the sequence,
in ascending order, of integers $i$ such that $0<i\le Q$ and $c_i=\tau(f_{i-1})$.
\end{itemize}

For consistency of notation in what follows, we set
$p(0)=p'(0)=q(0)=q'(0)=0$.

Thus each $b_i$, for $i=1,\dots,P$, can be written uniquely as
$b_{p(j)}$ or as $b_{p'(j)}$, and 
each $c_i$, for $i=1,\dots,Q$, can be written uniquely as
$c_{q(j)}$ or as $c_{q'(j)}$.

This notation allows us to give a more precise description of the
structure of the initial and terminal graphs of $\G$.    Specifically,
$I(\G)$ is constructed from the vertices $\{a,u,v\}$ by adding two edges
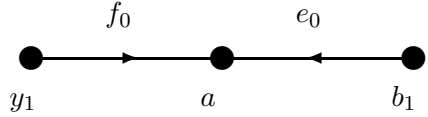
\eject

\vglue 10pt
\begin{center}
\setlength{\unitlength}{0.0125in}%
\begin{picture}(326,31)(155,495)
\thicklines
\put(240,520){\circle*{10}}
\put(320,520){\circle*{10}}
\put(400,520){\circle*{10}}
\put(240,520){\vector( 1, 0){ 45}}
\put(280,520){\line( 1, 0){ 40}}
\put(400,520){\vector(-1, 0){ 45}}
\put(320,520){\line( 1, 0){ 40}}
\put(230,500){\makebox(0,0)[lb]{\raisebox{0pt}[0pt][0pt]{$y_{1}$}}}
\put(390,500){\makebox(0,0)[lb]{\raisebox{0pt}[0pt][0pt]{$b_{1}$}}}
\put(310,500){\makebox(0,0)[lb]{\raisebox{0pt}[0pt][0pt]{$a$}}}
\put(270,535){\makebox(0,0)[lb]{\raisebox{0pt}[0pt][0pt]{$f_{0}$}}}
\put(350,535){\makebox(0,0)[lb]{\raisebox{0pt}[0pt][0pt]{$e_{0}$}}}
\end{picture}
\end{center}

together with directed chains

\vspace{20pt}
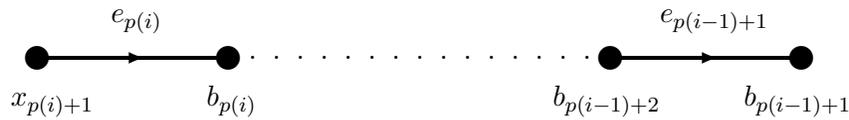
\begin{center}
\setlength{\unitlength}{0.0125in}%
\begin{picture}(326,31)(155,495)
\thicklines
\put(160,520){\circle*{10}}
\put(240,520){\circle*{10}}
\put(400,520){\circle*{10}}
\put(480,520){\circle*{10}}
\put(160,520){\vector( 1, 0){ 45}}
\put(200,520){\line( 1, 0){ 40}}
\multiput(240,520)(10.00000,0.00000){17}{\makebox(0.4444,0.6667){.}}
\put(400,520){\vector( 1, 0){ 45}}
\put(440,520){\line( 1, 0){ 40}}
\put(148,500){\makebox(0,0)[lb]{\raisebox{0pt}[0pt][0pt]{$x_{p(i)+1}$}}}
\put(230,500){\makebox(0,0)[lb]{\raisebox{0pt}[0pt][0pt]{$b_{p(i)}$}}}
\put(375,500){\makebox(0,0)[lb]{\raisebox{0pt}[0pt][0pt]{$b_{p(i-1)+2}$}}}
\put(455,500){\makebox(0,0)[lb]{\raisebox{0pt}[0pt][0pt]{$b_{p(i-1)+1}$}}}
\put(190,535){\makebox(0,0)[lb]{\raisebox{0pt}[0pt][0pt]{$e_{p(i)}$}}}
\put(420,535){\makebox(0,0)[lb]{\raisebox{0pt}[0pt][0pt]{$e_{p(i-1)+1}$}}}
\end{picture}
\end{center}

for each $i=1,\dots,s$, and

\vspace{20pt}
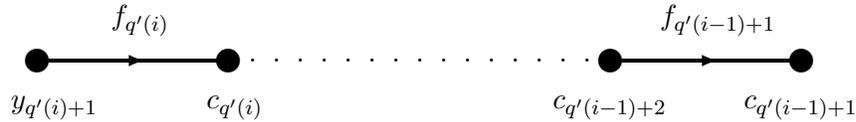
\begin{center}
\setlength{\unitlength}{0.0125in}%
\begin{picture}(326,31)(155,495)
\thicklines
\put(160,520){\circle*{10}}
\put(240,520){\circle*{10}}
\put(400,520){\circle*{10}}
\put(480,520){\circle*{10}}
\put(160,520){\vector( 1, 0){ 45}}
\put(200,520){\line( 1, 0){ 40}}
\multiput(240,520)(10.00000,0.00000){17}{\makebox(0.4444,0.6667){.}}
\put(400,520){\vector( 1, 0){ 45}}
\put(440,520){\line( 1, 0){ 40}}
\put(148,500){\makebox(0,0)[lb]{\raisebox{0pt}[0pt][0pt]{$y_{q'(i)+1}$}}}
\put(230,500){\makebox(0,0)[lb]{\raisebox{0pt}[0pt][0pt]{$c_{q'(i)}$}}}
\put(375,500){\makebox(0,0)[lb]{\raisebox{0pt}[0pt][0pt]{$c_{q'(i-1)+2}$}}}
\put(455,500){\makebox(0,0)[lb]{\raisebox{0pt}[0pt][0pt]{$c_{q'(i-1)+1}$}}}
\put(190,535){\makebox(0,0)[lb]{\raisebox{0pt}[0pt][0pt]{$f_{q'(i)}$}}}
\put(420,535){\makebox(0,0)[lb]{\raisebox{0pt}[0pt][0pt]{$f_{q'(i-1)+1}$}}}
\end{picture}
\end{center}

for each $i=1,\dots,t'$; and finally
single edges

\vspace{20pt}
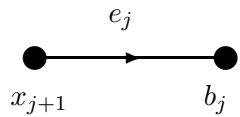
\begin{center}
\setlength{\unitlength}{0.0125in}%
\begin{picture}(326,31)(155,495)
\thicklines
\put(280,520){\circle*{10}}
\put(360,520){\circle*{10}}
\put(280,520){\vector( 1, 0){ 45}}
\put(320,520){\line( 1, 0){ 40}}
\put(269,500){\makebox(0,0)[lb]{\raisebox{0pt}[0pt][0pt]{$x_{j+1}$}}}
\put(350,500){\makebox(0,0)[lb]{\raisebox{0pt}[0pt][0pt]{$b_{j}$}}}
\put(310,535){\makebox(0,0)[lb]{\raisebox{0pt}[0pt][0pt]{$e_{j}$}}}
\end{picture}
\end{center}

for $p(s)<j\le P$ and

\vspace{20pt}
\begin{center}
\setlength{\unitlength}{0.0125in}%
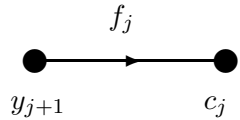
\begin{picture}(326,31)(155,495)
\thicklines
\put(280,520){\circle*{10}}
\put(360,520){\circle*{10}}
\put(280,520){\vector( 1, 0){ 45}}
\put(320,520){\line( 1, 0){ 40}}
\put(269,500){\makebox(0,0)[lb]{\raisebox{0pt}[0pt][0pt]{$y_{j+1}$}}}
\put(350,500){\makebox(0,0)[lb]{\raisebox{0pt}[0pt][0pt]{$c_{j}$}}}
\put(310,535){\makebox(0,0)[lb]{\raisebox{0pt}[0pt][0pt]{$f_{j}$}}}
\end{picture}
\end{center}

for $q'(t')<j\le Q$.

\bigskip
In the above diagrams $x_{P+1}$ and $y_{Q+1}$ (which have not been
defined) should be interpreted as $\iota(e_P)$ and $\iota(f_Q)$ respectively.
Note that at most one of these is equal to $a$.   (This happens if
and only if $a$ is the initial vertex of its incident edge in $\G$.)
All other $x_j$ and $y_j$ belong to $\{u,v\}$.

\bigskip
If $I(\G)$ contains a directed cycle, for example, then this cycle
must contain $a$.   From the above, we see that this can happen only if
$s=1$, $p(1)=P$, and $x_{P+1}=a$.

\bigskip
The structure of $T(\G)$ is entirely analogous, and
similar remarks apply.   We omit the details.

\bigskip
Now we are ready to construct a specific presentation for an HNN base for
$G=G(\G)$.   Recall that $G$ is given by a finite presentation
$${\cal P}(\G) = \< V(\G)~|~\iota(e)\lambda(e)=\lambda(e)\tau(e),~e\in E(\G)\>.$$
Since $\G$ is connected, we have $G^{ab}\cong\Z$, and the commutator
subgroup $G'$ is the normal closure in $G$ of the subgroup $B=B(\G)$
generated by the finite set $\{xy^{-1}~;~x,y\in V(\G)\}$.  A theorem of
Bieri and Strebel \cite{BS} says that $G$ is an HNN extension of $B$
with stable letter $t$ (which can be taken to be any element of $V(\G)$)
and associated subgroups $A_0=B\cap tBt^{-1}$ and $A_1=B\cap t^{-1}Bt$:
$$G=\<B,t~|~t^{-1}\alpha t=\phi(\alpha),~\alpha\in A_0\>,$$
where $\phi\co A_0\rightarrow A_1$ is the isomorphism induced by conjugation
by $t$.

Clearly $B$ is finitely generated.   It remains to prove that $B$ is
finitely presentable, and we do this by constructing an explicit 
set of defining relators.

Recall that our assumptions on $\G$
imply that each of $I(\G)$ and $T(\G)$ has precisely two connected
components, with the vertices $u,v$ belonging to separate components in
each case.

Let $F$ denote the subgroup of the free group on $V(\G)$ generated by
$$\{xy^{-1}~;~x,y\in V(\G)\}.$$   Then $F$ is free of rank 
$|V(\G)|-1=|E(\G)|$, and any basis for $F$ can be chosen as a finite
generating set for $B$.   Rather than fix a specific basis for $F$,
we proceed as follows.   Let $\bar K=\bar K(\G)$ be the maximal
abelian cover of the 2--complex $K=K(\G)$ associated to $\G$ (which
is the standard 2--complex model of the presentation ${\cal P}(\G)$).
Then since $K$ has a single 0--cell, we identify the 0--cells of
$\bar K$ with integers, via the isomorphism $H_1(K)\cong G^{ab}\cong\Z$.
The 1--cells of $\bar K$ with initial vertex $i\in\Z$ can be denoted
$w_i$, where $w\in V(\G)$, and each $w_i$ has terminal vertex $i+1\in\Z$.
Let $L$ be the 1--subcomplex of $\bar K$ with 0--cells $0,1$ and 1--cells
$\{w_0,~w\in V(\G)\}$.   Then $F$ is naturally identified with $\pi_1(L,0)$.

We also construct a graph $\hat L$ and an immersion $\pi\co \hat L\ra L$
as follows.   $V(\hat L)=\{0,1\}\times\{u,v\}$, $E(\hat L)=E(L)$,
$\iota(w_0)=(0,x)$ where $x\in\{u,v\}$ belongs to the same component
of $I(\G)$ as $w$, and $\tau(w_0)=(1,y)$ where $y\in\{u,v\}$ belongs
to the same component of $T(\G)$ as $w$.   The graph homomorphism
$\pi$ is defined to be the identity map on edges, and is defined on vertices by $\pi(i,u)=\pi(i,v)=i$, $i=0,1$.   It is not difficult to
see that $\hat L$ is connected.   Indeed, if the edge of $\G$ between
$u$ and $v$ has label $w$, then the edges $u,v,w$ of $\hat L$
form a spanning tree.
Since $\pi$ is bijective
on edges, it is an immersion, and hence injective on fundamental
groups.   Indeed, the fundamental group $\hat F$ of $\hat L$ embeds
as a free factor of $F=\pi_1(L)$ via $\pi_*$, as we can see by the
following construction: add an edge $X$ to $\hat L$ with $\iota(X)=(0,u)$
and $\tau(X)=(0,v)$, and an edge $Y$ with $\iota(Y)=(1,u)$,
$\tau(Y)=(1,v)$, to form a larger graph $\tilde L$.   The immersion
$\pi\co \hat L\ra L$ extends to a homotopy equivalence $\pi\co \tilde L\ra L$
that shrinks the edge $X$ to the vertex $0$, and the edge $Y$ to the
vertex $1$.   Hence we have 
$$F=\pi_1(L)\cong\pi_1(\tilde L)=\pi_1(\hat L)*\<X,Y\>.$$

Since the map $\pi\co \hat L\ra L$ is bijective on edges, any path in $L$
which lifts to a path in $\hat L$ does so uniquely.   Given a closed
path $C$ in $L$ that lifts to a closed path $\hat C$ in $\hat L$, we
define two related paths in $L$, namely the {\it forward derivative}
$\d_+C$ of $C$ and the {\it backward derivative} $\d_-C$ of $C$, as
follows.  For $\d_+C$ we first fix a maximal subforest $\Phi_I$
of $I(\G)$.   Next, we cyclically permute $\hat C$ so that it
begins and ends at one of the vertices $(1,u)$ or $(1,v)$.   Hence $\hat C$
is a concatenation of length two subpaths of the form $x^{-1}y$, where
$x,y\in E(\hat L)=V(\G)$ belong to the same component of $I(\G)$.
The next step is to replace each such subword $x^{-1}y$ by the product
$$(x^{-1}z_0)(z_0^{-1}z_1)\dots(z_m^{-1}y),$$
where $(x,z_0,z_1,\dots,z_m,y)$ is the geodesic from $x$ to $y$ in $\Phi_I$.   We now have a concatenation of length 2 subwords of the
form $x^{-1}y$ where $x$ and $y$ are joined by an edge in $\Phi_I$.  
This edge corresponds to an edge of $\G$, and hence to a defining
relation in ${\cal P}(\G)$ that can be written
$$x^{-1}y=gh^{-1}$$
for some $g,h\in V(\G)$.  The final step is to replace each such 
word $x^{-1}y$ by the corresponding word $gh^{-1}$.   The result is
a closed path $\d_+C$ in $L$.

\begin{remarks}
\begin{enumerate}
\item $\d_+C$ depends on the choice of maximal forest $\Phi_I$, and
then is well-defined only up to cyclic permutation.
\item If $C'$ is a cyclic permutation of $C$, then $C'$ also lifts to
a closed path in $\hat L$, so $\d_+C'$ is defined.   It is equal to
(a cyclic permutation of) $\d_+C$.
\item The definition of $\d_+C$ does not depend on $C$ being (cyclically)
reduced.   Indeed the insertion into $C$ of a cancelling pair $xx^{-1}$
may alter $\d_+C$.   However, the insertion of a cancelling pair
$x^{-1}x$ will {\it not} alter $\d_+C$.
\item $C$ and $\d_+C$ are (freely) homotopic in $\bar K$ (since the 
last part of the construction involves replacing a path 
$x^{-1}y$ by a homotopic path $gh^{-1}$).   In particular, if $C$ is
nullhomotopic in $\bar K$, then so is $\d_+C$.
\item The unique lift of $\d_+C$ in $\tilde L$ does not contain the
edge $Y$.
\end{enumerate}
\end{remarks}

The backward derivative $\d_-C$ is defined similarly.   This time we
fix a maximal forest $\Phi_T$ of $T(\G)$, and choose a cyclic
permutation of $\hat C$ beginning at $(0,u)$ or $(0,v)$, split $\hat C$
into subpaths of the form $xy^{-1}$ with $x,y$ in the same component of $T(\G)$,
and then use relations of ${\cal P}$ corresponding to edges of $\Phi_T$
to transform $\hat C$.   Remarks analogous to the above hold
also for $\d_-C$.

\bigskip
Now consider the unique cycle in $T(\G)$.   If $z_0,\dots,z_m$
are the vertices of this cycle in cyclic order, define $\hat R_0$
to be the nullhomotopic path 
$$(z_mz_0^{-1})(z_0z_1^{-1})\dots(z_{m-1}z_m^{-1})$$
in $\hat L$ and $R_0=\pi(\hat R_0)$ the corresponding nullhomotopic
path in $L$.   Now define $R_1=\d_-R_0$. If $R_1$ lifts to $\hat L$
then define $R_2=\d_-R_1$, and so on.   In this way we obtain either
an infinite sequence $R_1,R_2,\dots$ of paths in $L$, or a finite
sequence $R_1,\dots,R_M$ such that $R_M$ does not lift to $\hat L$.

In a similar way, the unique cycle in $I(\G)$ determines a nullhomotopic
closed path $S_0$ in $L$ that lifts to $\hat L$, so a sequence
$S_1,\dots$ of closed paths in $L$ (finite or infinite), such that
$S_i=\d_+S_{i-1}$ for each $i\ge 1$, and if the sequence is finite
with final term $S_N$ then $S_N$ does not lift to $\hat L$.

\begin{lemma}\sl
The paths $R_i$ and $S_j$ are all nullhomotopic in $\bar K$.
\end{lemma}

\begin{proof}
This follows by induction and Remark (iv) above, since $R_0$ and $S_0$
are nullhomotopic.
\end{proof}

Now suppose that the sequence $\{R_i\}$ contains at least $m$ terms.
We construct a 2--complex $L_m$ as follows.   The 1--skeleton of $L_m$
is the subcomplex of $\bar K$ consisting of $L$, together with the
0--cells $2,\dots,m+1$ and the 1--cells $u_1,v_1,\dots,u_m,v_m$.
Then $L_m$ has precisely $m$ 2--cells attached to $L$ using the
paths $R_1,\dots,R_m$.   We also consider the full subcomplex
$\bar K_m$ of $\bar K$ on the set $\{0,1,\dots,m+1\}$ of 0--cells.

\begin{lemma}\sl
The 2--complexes $L_m$ and $\bar K_m$ are homotopy equivalent.
\end{lemma}

\begin{proof}
We argue by induction on $m$, there being nothing to prove in the
case $m=0$.   Let $\gamma$ denote the covering transformation of
$\bar K$ that sends a 0--cell $n\in\Z$ to $n+1$. 
Note that the link of the 0--cell $m+1$ in $\bar K_m$ is naturally
identifiable with the graph $T(\G)$.   Let $d$ be the unique
edge in $E(\G)=E(T(\G))$ that does not belong to the maximal forest
$\Phi_T\subset T(\G)$.   Then $d$ is contained in the unique cycle in
$T(\G)$, so $R_0$ has a subword $xy^{-1}$, where $x,y$ are the endpoints
of $d$ in $T(\G)$.
Corresponding to $d$ is a relator $xy^{-1}h^{-1}g$ in $\cal P$,
which lifts to a 2--cell $\alpha$ with boundary path
$x_my_m^{-1}h_{m-1}^{-1}g_{m-1}$ in $\bar K_m$.   Modulo the other
2--cells of $\bar K_m$, the boundary path of $\alpha$ is homotopic
to $\gamma^m(R_0)^{-1}\cdot\gamma^{m-1}(R_1)$.   Since $R_0$ is nullhomotopic
in the 1--skeleton of $\bar K$, this is in fact homotopic to
$\gamma^{m-1}(R_1)$.   This in turn is homotopic (in $\bar K_{m-1}$) to
$\gamma^{m-2}(R_2)$, etc.   Repeating this argument, we see that
the boundary path of $\alpha$ is homotopic in $\bar K_m\backslash\alpha$
to $R_m$.   A simple homotopy move allows us to replace $\alpha$
by a 2--cell whose boundary path is $R_m$.

The link of $m+1$ in the resulting 2--complex $K'$ is then isomorphic to
$T(\G)\backslash d=\Phi_T$.   Since $\Phi_T$ is a forest with two
components (one containing $u$ and the other containing $v$),
it collapses to the graph with no edges and vertex set $\{u,v\}$.
Each move in this collapsing process (removing a vertex and an edge
from the graph) can be mirrored by a collapse in the 2--complex $K'$
(removing a 1--cell and a 2--cell that are incident at the 0--cell $m+1$).
After performing all these collapsing moves, we are left with a 2--complex
$K''$, simple homotopy equivalent to $\bar K_m$.   By inspection,
$K''$ is formed from $\bar K_{m-1}$ by adding a 2--cell with boundary
path $R_m$, a 0--cell $m+1$, and two 1--cells $u_m,v_m$, each joining
$m$ to $m+1$.

By inductive hypothesis, $\bar K_{m-1}$ is homotopy equivalent to
$L_{m-1}$, so $\bar K_m$ is homotopy equivalent to the 2--complex
obtained from $L_{m-1}$ by adding a  2--cell with boundary
path $R_m$, a 0--cell $m+1$, and two 1--cells $u_m,v_m$, each joining
$m$ to $m+1$.   But this 2--complex is precisely $L_m$, and the proof is complete.
\end{proof}

\begin{remark}
An analogous result holds for the $S_j$.   We omit the details, but 
will use this result implicitly in what follows.
\end{remark}

\begin{cor}\sl
If $R_1,\dots,R_m$ and $S_1,\dots,S_n$ are all defined, then
$m+n<|V(\G)|$.
\end{cor}

\begin{proof}
By the Lemma and its analogue for the $S_j$, $\bar K_m$ is homotopy
equivalent to a 2--complex formed from $L$ by attaching $m$ 2--cells and
then wedging on $m$ circles; and $\gamma^{-n}(\bar K_n)$ is homotopy
equivalent to a complex obtained from $L$ by adding $n$ 2--cells
and then wedging on $n$ circles.   Since
$\gamma^{-n}(\bar K_{m+n})=\gamma^{-n}(\bar K_n)\cup\bar K_m$,
with $\gamma^{-n}(\bar K_n)\cap\bar K_m=\bar K_1=L$, it follows
that $\gamma^{-n}(\bar K_{m+n})$ is homotopy equivalent to
a complex formed from $L$ by adding $m+n$
2--cells and then wedging on $m+n$ circles.   Hence 
$\beta_1(\bar K_{m+n})\ge m+n$.
Now $H_2(K)=0$, and $\bar K$ is a $\Z$--cover of $K$, so $H_2(\bar K)=0$
by \cite{A}, Proposition 1.  Hence also $H_2(K')=0$ for any subcomplex
$K'\subseteq K$.   In particular
 $H_2(\bar K_{m+n})=0=H_2(L)$.
Since also $H_0(\bar K_{m+n})=\Z=H_0(L)$ and 
 $\chi(\bar K_{m+n})=\chi(L)=2-|V(\G)|$,
it follows that 
$$m+n\le\beta_1(\bar K_{m+n})=\beta_1(L)=|V(\G)|-1.$$

\vspace{-2\baselineskip}\end{proof}

\begin{cor}\sl
Each of the sequences $\{R_i\}$ and $\{S_j\}$ are finite, and if the
final terms are $R_M$ and $S_N$ respectively then $M+N<|V(\G)|$.
\end{cor}

We claim  that the finite sequences $\{R_i\}$ and $\{S_j\}$ form
a full set of defining relators for the HNN base $B$ of $G$, which
completes the proof of our Theorem \ref{main}.   In order to prove this
claim, we need to derive some further information about the structure
of the words $R_i$ and $S_j$.

\begin{remark}
The definitions of $R_i$ and $S_i$ depend, {\it a priori}, on specific
choices for the maximal forests $\Phi_T$ and $\Phi_I$ respectively.
Suppose we were to choose a different maximal tree $\Phi_I'$ in $I(\G)$,
for example.   Then geodesics in $\Phi_I$ and $\Phi_I'$ would differ
at most by the unique cycle in $I(\G)$.  It follows from this that
the resulting definitions of $\d_+C$, for any closed path $C$ in $L$
that lifts to $\hat L$, are equal modulo the normal closure of $S_1$.
An easy induction shows that, for any $i$, the definitions of $S_i$ resulting from different
choices of $\Phi_I$ are equal modulo the normal closure of $\{S_1,\dots,S_{i-1}\}$.
Hence our set of defining relators does not depend in an essential
way upon the choices of maximal forests $\Phi_I$ and $\Phi_T$.
\end{remark}

\section{Structure of the relations}\label{struct}

In this section we examine the structure of the proposed defining
relators $R_i$
and $S_i$ of the HNN base $B$ for $G$.   Recall that
each of $R_i$ and $S_i$ is a closed path in the 2--complex $L$, and that
we have a homotopy equivalence $\pi\co \tilde L\ra L$, which restricts
to an edge-bijective
graph immersion on $\hat L=\tilde L\backslash\{X,Y\}$ and shrinks
each of the 1--cells $X,Y$ to a point.   Let $\tilde C$ denote the unique
(up to cyclic permutation) cyclically reduced closed path in $\tilde L$ 
that maps to a given cyclically reduced closed path $C$ in $L$.   Then
$C$ lifts to $\hat L$ if and only if $\tilde C$ is a path in $\hat L$, 
in which case $\tilde C$ is the unique lift.   By definition, each
$R_i$ (resp $S_i$) is defined if and only if $R_{i-1}$ (resp $S_{i-1}$)
lifts to $\hat L$.   Hence $\tilde R_i$ is a path in $\hat L$ for
$1\le i\le M-1$, and $S_i$ is a path in $\hat L$ for $1\le i\le N-1$.
Moreover, the path $\tilde R_M$ involves $Y$ but not $X$, while the path
$\tilde S_N$ involves $X$ but not $Y$.

\bigskip For any group $A$ and letter $Z$, we say that a word
$w\in A*\<Z\>$ is {\it positive} (resp {\it negative}) in $Z$ if
only positive (resp negative) powers of $Z$ occur in $w$.
We say that $w$ is {\it strictly positive} (resp {\it strictly
negative}) if in addition at least one positive (resp negative)
power of $Z$ does occur in $w$, in other words $w\not\in A$.

\bigskip
We will concentrate our attention on the relators $S_i$.   The analysis
of the $R_i$ is entirely analogous.

We first treat the case where $I(\G)$ contains a directed cycle $C$.

\begin{theorem}\label{struct0}\sl
Suppose that the unique cycle $C$ in $I(\G)$ is directed. Then:
\begin{itemize}
\item $N=1$;
\item $\tilde S_1$ is either strictly positive or strictly negative in $X$;
\item $S_1$ involves each of $a,b_1,\dots,b_P$ exactly once, and no $c_j$;
\item each of $a,b_1,\dots,b_P$ is an extremal source in $\G$.
\end{itemize}
\end{theorem}

\begin{proof}
The vertex $a$ is contained in $C$, by Lemma \ref{I}, (v). 
Since $\iota(f_0)\in\{u,v\}$, $f_0$ is not an edge of $C$, so the edge
of $C$ coming into $a$ is $e_0$.   Hence $b_1=\iota(e_0)$ is a vertex
of $C$, and since $e_1$ is the only edge with $\lambda(e_1)=b_1$,
it is also an edge of $C$, and so on.   Hence each of $b_1,\dots,b_P$
are vertices of $C$, $\iota(e_P)=a$,
and the edges of $C$ are precisely $e_P,\dots,e_0$
(in the order of the orientation of $C$).
Each of the vertices of $C$ is
extremal in $\G$, and since it is the initial vertex of an edge of $I(\G)$
it is also the initial vertex of an edge of $\G$, ie a source in $\G$.
Moreover 
$$S_0=(a^{-1}b_P)(b_P^{-1}b_{P-1})\dots(b_1^{-1}a),$$
so
$$S_1=\d_+S_0=(b_P\tau(e_P)^{-1})(b_{P-1}x_P^{-1})\dots(b_1x_2^{-1})(ax_1^{-1}),$$
where each $x_i\in\{u,v\}$.

Suppose that $S_1$ lifts to $\hat L$.   Then $\tau(e_P)$ belongs to the same
component of $I(\G)$ as $b_{P-1}$, $x_P$ to the same component as $b_{P-2}$,
and so on.   Since $a,b_1,\dots,b_P$ all belong to the same component
of $I(\G)$, it follows that the $x_i$ also all belong to the same
component.   But $u$ and $v$ belong to different components of $I(\G)$,
and so the $x_i$ are all equal, which contradicts Lemma \ref{key}.

Hence $S_1$ does not lift to $\hat L$, and so $N=1$.   Moreover, by
the above argument, some of the $x_i$ belong to the opposite component 
of $I(\G)$ from $a$.   If $a,u$ belong to the same component of $I(\G)$,
this means that some of the $x_i$ are equal to $v$.
Then $\tilde S_1$ is formed from $S_1$ by replacing each occurrence
of $v^{-1}$ by $v^{-1}X^{-1}$, and so $\tilde S_1$ is strictly
negative in $X$.   Similarly, if $a,v$ belong to the same component of
$I(\G)$, then $\tilde S_1$ is strictly positive in $X$.
\end{proof}

For the rest of the section, we can assume that the cycle $C$ is not
directed.   Then $y_1=\iota(f_0)=\iota(e_{p(1)})\in\{u,v\}$.
We may assume that $y_1=u$.   Then $C$ has the form

\begin{center}
\setlength{\unitlength}{0.0125in}%
\begin{picture}(350,127)(140,490)
\thicklines
\put(160,520){\circle*{10}}
\put(240,520){\circle*{10}}
\put(400,520){\circle*{10}}
\put(480,520){\circle*{10}}
\put(160,600){\circle*{10}}
\put(480,600){\circle*{10}}
\put(160,520){\vector( 1, 0){ 45}}
\put(200,520){\line( 1, 0){ 40}}
\multiput(240,520)(10.00,0.00000){17}{\makebox(0.4444,0.6667){.}}
\put(400,520){\vector( 1, 0){ 45}}
\put(440,520){\line( 1, 0){ 40}}
\put(160,600){\vector( 1, 0){165}}
\put(480,520){\vector( 0, 1){ 45}}
\put(160,600){\vector( 0,-1){ 45}}
\put(160,560){\line( 0,-1){ 40}}
\put(320,600){\line( 1, 0){160}}
\put(480,560){\line( 0, 1){ 40}}
\put(485,610){\makebox(0,0)[lb]{\raisebox{0pt}[0pt][0pt]{$a$}}}
\put(150,490){\makebox(0,0)[lb]{\raisebox{0pt}[0pt][0pt]{$b_{p(1)}$}}}
\put(225,490){\makebox(0,0)[lb]{\raisebox{0pt}[0pt][0pt]{$b_{p(1)-1}$}}}
\put(395,490){\makebox(0,0)[lb]{\raisebox{0pt}[0pt][0pt]{$b_2$}}}
\put(475,490){\makebox(0,0)[lb]{\raisebox{0pt}[0pt][0pt]{$b_1$}}}
\put(140,610){\makebox(0,0)[lb]{\raisebox{0pt}[0pt][0pt]{$u$}}}
\end{picture}
\end{center}

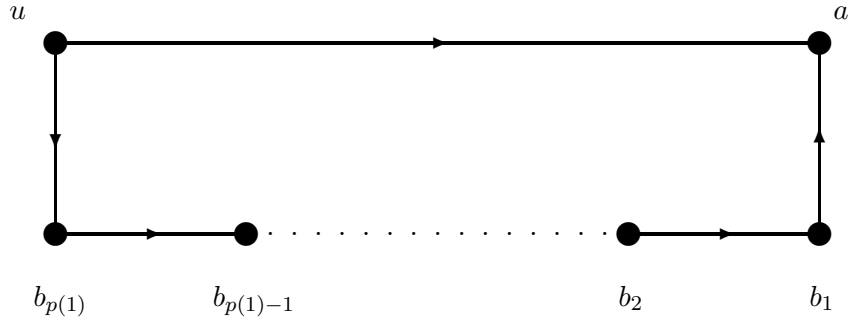
\centerline{\small Figure 1}

\bigskip
For the purpose of defining forward derivatives, and hence the $S_i$,
we fix $\Phi_I$ to be the maximal subforest of $I(\G)$ obtained by
removing the edge $f_0$ (the edge joining $u$ to $a$ in $C$).

For $k\le{\rm min}(s,t'+1)$, let $I_k(\G)$ denote the subgraph of $\Phi_I$
consisting of the edges $\{e_i,0\le i\le p(k)\}$ and $\{f_i,1\le i\le q'(k-1)\}$,
together with all their incident vertices.
Note that $I_k$ contains no more than two components, one contained in
each component of $\Phi_I$.   Hence whenever two vertices of $I_k$
belong to the same component of $\Phi_I$, then the geodesic between them
is also contained in $I_k$.

\begin{theorem}\label{struct1}\sl
Suppose that the cycle in $I(\G)$ has the form shown in Figure 1.
Then:
\begin{enumerate}
\item
Each $S_i$ can be written, up to cyclic permutation, in the form
$aU_ia^{-1}V_i$, where $U_i$ is a word in 
$$\{a,u,v,c_1,\dots,c_{q'(i-1)+1}\};$$ and $V_i$ is a word in $$\{a,u,v,b_1,\dots,b_{p(i)+1}\}.$$
\item
If $p(i)<P$, then $V_i$ contains a single occurrence of
$b_{p(i)+1}$ and does not contain $a$.
\item
If $q'(i-1)<Q$, then $U_i$ contains a single occurrence of  $c_{q'(i-1)+1}$
 and does not contain $a$.
\item
Every letter occurring in $S_i$, other than $b_{p(i)+1}$ and $c_{q'(i-1)+1}$,
 is a vertex of the subgraph
$I_i\subseteq I(\G)$.
\item
If $p(i)=P$ or $q'(i-1)=Q$ then $i=N$.
\end{enumerate}
\end{theorem}

\begin{proof}
We prove this by induction on $i$, the initial case being when $i=1$.
We have 
$$S_0=(u^{-1}a)(a^{-1}b_1)(b_1^{-1}b_2)\dots(b_{p(1)}^{-1}u),$$
so
$$S_1=\d_+S_0=(ac_1^{-1})(x_1a^{-1})(x_2b_1^{-1})\dots(x_{p(1)}b_{p(1)-1}^{-1})(b_{p(1)+1}b_{p(1)}^{-1})$$
(if $p(1)<P$).   The vertices $a,u,b_1,\dots,b_{p(1)}$ are contained in $I_1$, but not $c_1$, $b_{p(1)+1}$.   The first four statements of the result (for $i=1$) follow, setting $U_1=c_1^{-1}x_1$ and
$$V_1=(x_2b_1^{-1})\dots(x_{p(1)}b_{p(1)-1}^{-1})(b_{p(1)+1}b_{p(1)}^{-1}).$$

For the last statement, certainly $Q>0=q'(0)$.   Suppose that $p(1)=P$ and
$i<N$.   Then 
$$S_1=(ac_1^{-1})(x_1a^{-1})(x_2b_1^{-1})\dots(x_Pb_{P-1}^{-1})(\tau(e_P)b_P^{-1})$$
lifts to $\hat L$, so each of $x_2,\dots,x_P$
belongs to the same component of $I(\G)$ as $a,b_1,\dots,b_{P-1}$,
in other words $x_2=\dots=x_P=u$. By Lemma \ref{key} we have
$x_1=v$ and $e_P$ incident with $v$.  But $\iota(e_P)=u$ so $\tau(e_P)=v$,
which does not belong to the same component of $I(\G)$ as $b_{P-1}$.
It follows that $S_1$ does not, after all, lift to $\hat L$, a contradiction.

This completes the proof of the initial
case of the induction.

\bigskip
Now assume inductively that $i>1$ and the result is true for $i-1$.
In particular, $i-1<N$, so $p(i-1)<P$ and $q'(i-2)<Q$.
Hence $U_{i-1}$ contains a single occurrence of $c_{q'(i-2)+1}$,
$V_{i-1}$ contains a single occurrence of $b_{p(i-1)+1}$, and every
other letter occurring in $S_{i-1}$ is a vertex of the subgraph
$I_{i-1}$ of $I(\G)$.   Consider the construction of $S_i=\d_+S_{i-1}$
from $S_{i-1}$.   We first write a suitable cyclic permutation of
$S_{i-1}$ as a product of length two subwords of the form $g^{-1}h$.
For all but two of these subwords, both $g$ and $h$ are vertices of 
$I_{i-1}$.  (There are precisely two exceptions, since the occurrences of 
$b_{p(i-1)+1}$ and $c_{q'(i-2)+1}$ in $S_{i-1}$ are separated at least
by an occurrence of $a^{\pm 1}$.)

Suppose first that $g,h$ are vertices of $I_{i-1}$.   The next step
is to replace $g^{-1}h$ by the product 
$$(g^{-1}z_1)(z_1^{-1}z_2)\dots(z_t^{-1}h)$$
where $g,z_1,z_2,\dots,z_t,h$ are the vertices on the geodesic from $g$
to $h$ in $\Phi_I$.   This geodesic is contained in $I_{i-1}$,
so each bracketed term here is $(\iota(e)^{-1}\lambda(e))^{\pm 1}$ for some edge $e$ of $I_{i-1}$.
The final step is to replace this by $(\lambda(e)\tau(e)^{-1})^{\pm 1}$.
Note that $\tau(e)$ is a vertex of $I_i$, and $\tau(e)\ne a$.
Also, none of the intermediate vertices $z_i$ in the geodesic is equal
to $a$, since $a$ is an extremal vertex of $\Phi_I$.
Note that, if $g^{-1}h$ is a subword of $U_{i-1}$, then all letters
in the resulting subword of $S_i$ come from $\{u,v,c_1,\dots,c_{q'(i-1)}\}$,
while if it is a subword of $a^{-1}V_{i-1}a$ then all letters come from
$\{a,u,v,b_1,\dots,b_{p(i)}\}$.

A similar argument holds if, say $g=b_{p(i-1)+1}$.   Here, however,
the geodesic from $g$ to $h$ is not contained in $I_{i-1}$.   It is
the union of the geodesic from $b_{p(i-1)+1}$ to $z$ in $I_i$,
where $z\in\{u,v\}$, with the geodesic (in $I_{i-1}$) from $z$ to $h$.
Edges in $I_{i-1}$ give rise to length 2 subwords of $S_i$ consisting
of letters which are vertices in $I_i$.   The same is true for an
edge $e_j$ from $b_j$ to $b_{j+1}$, for $p(i-1)<j<p(i)$.
(The corresponding word is $x_jb_j^{-1}$.)   Finally, the edge $e_{p(i)}$
(from $b_{p(i)}$ to $z$) contributes a subword $\tau(e_{p(i)})b_{p(i)}^{-1}$.
If $p(i)<P$ then  $\tau(e_{p(i)})=b_{p(i)+1}$; otherwise $\tau(e_{p(i)})\in\{a,u,v\}$.

The analysis if $h=b_{p(i-1)+1}$, or if one of $g,h$ is $c_{q'(p-2)+1}$
is similar to the above.

\bigskip
Each of the two subwords $g^{-1}h$ of $S_{i-1}$ that contain the letter
$a$ gives rise to a subword of $S_i$ containing an occurrence of $a$
with the same exponent.   If $g=a$ then the subword
begins $(x_1a^{-1})\dots$, while if $h=a$ then the subword ends
$\dots(ax_1^{-1})$.   If $p(i)<P$ and $q'(i-1)<Q$ then this will be the only 
occurrence of $a$ in this subword of $S_i$.

\bigskip
Statements (i)--(iv) follow.

\bigskip
To prove (v), suppose for example that $i<N$ and $p(i)=P$.   
Another induction
on $i$ shows that $x_2=\dots=x_P=u$.  An argument
similar to that given above in the initial case of the induction
again gives rise to a contradiction: by Lemma \ref{key}, $\tau(e_P)=v$,
which does not belong to the same component of $I(\G)$ as $b_{P-1}$,
so $S_i$ does not lift to $\hat L$ and $i=N$.

If $i<N$ and $q'(i-1)=Q$ then a similar argument applies.   Here
we can show that $y_1=\dots=y_Q=x_1\in\{u,v\}$, which contradicts
Lemma \ref{key}.
\end{proof}

This result contains all the necessary information about $S_i$ if $i<N$.
We now need to investigate further the structure of $\tilde S_N$, particularly
as regards occurrences of $X$.
Note that, up to cyclic permutation, we have
$\tilde S_N=a\tilde U_Na^{-1}\tilde V_N$, by Theorem \ref{struct1} (i).

\begin{lemma}\label{SNX}\sl
Each of $\tilde U_N$, $\tilde V_N$ is either positive or negative
in $X$.
\end{lemma}

\begin{proof}
As indicated in the proof of Theorem \ref{struct1}, all of $V_N$, except for the
part arising from the geodesic $\gamma$
from $b_{p(N-1)+1}$ to $u$, consists
of letters which are vertices in $I_{N-1}$.
All of these vertices are in the same component of $I(\G)$
as $u$.   The part of $V_N$ arising from $\gamma$ is 
$$[(x_{p(N-1)+2}b_{p(N-1)+1}^{-1})\dots(x_{p(N)}b_{p(N)-1}^{-1})(\tau(e_{p(N)})b_{p(N)}^{-1})]^{\pm 1},$$
or, if $\gamma$ passes through $a$ (ie if $\iota(e_{p(N)})=a$):
$$[(x_{p(N-1)+2}b_{p(N-1)+1}^{-1})\dots(\tau(e_{p(N)})b_{p(N)}^{-1})(x_1a^{-1})\dots(b_{p(1)+1}b_{p(1)}^{-1})]^{\pm 1}.$$
The expression in square brackets is a product of terms $gh^{-1}$ with $h$ in the same
component of $I(\G)$ as $u$.   To lift to $\tilde L$, we replace $h^{-1}g$
by $h^{-1}Xg$ whenever $g$ belongs to the same component of $I(\G)$
as $v$ and $h$ to the same component as $u$, and by $h^{-1}X^{-1}g$
if $g$ belongs to the same component as $u$ and $h$ to the same component
as $v$.   Hence $\tilde V_N$ is either positive or negative in $X$

A similar argument applies to $\tilde U_N$, replacing $u$ by $x_1$ in the
above.  
\end{proof}

We will also need to investigate possible occurrences of $a$ in
$S_N$ other than those indicated in Theorem \ref{struct1}.   

\begin{lemma}\label{SNa}\sl
The words $\tilde U_N$ and $\tilde V_N$ contain in total at most one occurrence of $a$.
\end{lemma}

\begin{proof}
From the discussion in the proof of Lemma \ref{SNX}, the word $V_N$
(and hence also $\tilde V_N$)
contains a single occurrence of $a$ if $e_{p(N)}$ is incident with $a$ in $\G$, and 
no occurrence of $a$ otherwise.   Similarly $U_N$ 
(and hence also $\tilde U_N$) contains a single occurrence of $a$ if $f_{q'(N-1)}$ is incident with $a$ in $\G$, and no occurrence of $a$
otherwise.   The result now follows from the fact that $a$ is extremal in $\G$.
\end{proof}

\section{Completion of the proof}\label{kill}

Define 
$$G_0=\pi_1(\hat L)/\{R_1,\dots,R_{M-1},S_1,\dots,S_{N-1}\},$$
$$G_+=(G_0*\<X\>)/\{\tilde S_N\},$$
$$G_-=(G_0*\<Y\>)/\{\tilde R_M\},$$
and
$$G_1=(G_0*\<X,Y\>)/\{\tilde R_M,\tilde S_N\}\cong(\pi_1(L))/\{R_1,\dots,R_M,S_1,\dots,S_N\}.$$

\begin{lemma}\label{free}\sl
The group $G_0$ is free. 
\end{lemma}

\begin{proof}
By Theorems \ref{struct0} and \ref{struct1}, and the analogous results for the $R_i$,
the set of $M+N-2$ distinct numbers 
${\cal B}=\{p(1)+1,\dots,p(N-1)+1,p'(0)+1,\dots,p'(M-2)+1\}$
has the property that each $j\in{\cal B}$ is the greatest index of a
$b$--letter occurring in a unique relator $R_i$ or $S_i$, and moreover
that relator contains precisely one occurrence of $b_j$.   

It follows that the 1--complex $L'$ obtained from $\hat L$ by removing the
1--cells $b_j,~j\in{\cal B}$ is connected, with fundamental group
isomorphic to $G_0$. 
\end{proof}

\begin{lemma}\label{injectivity}\sl
The natural maps $G_0\ra G_+$ and $G_0\ra G_-$ are injective.
\end{lemma}

\begin{proof}
We show that the map $G_0\ra G_+$ is injective.   The proof
of injectivity of $G_0\ra G_-$ is entirely analogous.
Since $G_0$ is a free group and $G_+$ is a one-relator group
$G_+=(G_0*\< X\>)/\{\tilde S_N\}$, we need only show that $\tilde S_N$,
regarded as a word in $(G_0*\< X\>)$, genuinely involves $X$.
The result then follows from the Freiheitssatz for one-relator
groups \cite{Mag}.

Consider the various possibilities for the structure of $\tilde S_N$.
If the initial graph $I(\G)$ contains a directed cycle, then $N=1$
and $\tilde S_1$ is a strictly positive (or strictly negative)
word in $X$, by Theorem \ref{struct0}.   Thus $\tilde S_1$,
regarded as a word in the free product $G_0*\< X\>$, is also
strictly positive (or strictly negative)
in $X$, and so genuinely involves $X$.

\bigskip
Suppose then that $I(\G)$ does not contain a directed cycle.
By Theorem \ref{struct1} (i) and Corollary \ref{SNX}
we have (up to cyclic permutation)
$\tilde S_N=a\tilde U_Na^{-1}\tilde V_N$, with each of $\tilde U_N$ and $\tilde V_N$ being either
positive or negative in $X$.   We also have $\tilde S_N$ definitely involving
$X$, since otherwise $S_N$ would lift to $\hat L$.

If $X$ occurs in $\tilde S_N$ with nonzero exponent-sum, then occurrences of
$X$ survive modulo the relators $R_1,\dots,R_{M-1},S_1,\dots,S_{N-1}$, so we may
assume that $X$ appears with exponent-sum zero.   Thus one of $\tilde U_N$, $\tilde V_N$
is strictly positive, and the other is strictly negative, with 
precisely the same number of occurrences of $X^{\pm 1}$.
We may rewrite $\tilde S_N$ (again, up to cyclic permutation) as
$$\tilde S_N=XA_1X\dots A_tXW_1X^{-1}B_tX^{-1}\dots B_1X^{-1}W_2$$
for some $t\ge 0$ and words $A_i,B_i$ and $W_1,W_2$ that do not
involve $X$.   If we can show that neither $W_1$ nor $W_2$ is equal to
the identity element in $G_0$, then it will follow that the above
expression for $\tilde S_N$ does not allow for cancellation of $X$--symbols,
when reducing modulo the relators of $G_0$.   The result will follow.

Now $a$ occurs with exponent-sum zero in each of the relators
$R_1,\dots,R_{M-1}$ and $S_1,\dots,S_{N-1}$ of the group $G_0$,
by Theorem \ref{struct1}.   If neither
$U_N$ nor $V_N$ contains the letter $a$, then each of $W_1$, $W_2$
contains precisely one occurrence of $a$, and so has infinite order in $G_0$.
In particular, they are nontrivial in $G_0$, as required.

This reduces us to the case where one of $U_N$, $V_N$ involves the letter 
$a$.  By Corollary \ref{SNa} we know that this can happen for only one of
$U_N$, $V_N$.

\bigskip
First
suppose that $a$ occurs in $U_N$.   Then $q'(N-1)=Q$
(and so also $N>1$).   As in the proof
of Corollary \ref{SNX}, the part of $U_N$ that gives rise to occurrences of $X$
comes from the geodesic $\delta$ in $\Phi_I$ from $c_{q'(N-2)+1}$ to $x_1$.
The relevant subword of $U_N$ has the form:
$$[(y_{q'(N-2)+2}c_{q'(N-2)+1}^{-1})\dots(y_{Q}c_{Q-1}^{-1})(\tau(f_Q)c_Q^{-1})]^{\pm 1},$$
or, if $\delta$ passes through $a$:
$$[(y_{q'(N-2)+2}c_{q'(N-2)+1}^{-1})\dots(\tau(f_Q)c_Q^{-1})(x_1a^{-1})\dots(b_{p(1)+1}b_{p(1)}^{-1})]^{\pm 1}.$$
The occurrences of $X$ in $\tilde U_N$
correspond to those $y_j$, $j\ge q'(N-2)+2$
that are not equal to $x_1$, and also from $\tau(f_Q)$ if this is not
in the same component of $I(\G)$ as $x_1$.   In the case where $\delta$ passes through 
$a$, we see that, in $\tilde S_N=a\tilde U_Na^{-1}\tilde V_N$ the $a$--letters that occur in 
the same $W_i$ have the same exponent, and hence the $W_i$ are both
nontrivial in $G_0$, as required.
In the other case, $\tau(f_Q)=a$ and the unique occurrence of $c_Q$ in $\tilde V_N$
lies on the same side of all the $X$--letters as the unique occurrence of
$a$.  Hence
$c_Q$ occurs (precisely once) in the same $W_i$ that contains two $a$--letters.   To prove that this $W_i$ is nontrivial
in $G_0$, it suffices to show that $c_Q$ does not occur in any of the
relators $R_1,\dots,R_{M-1}$ or $S_1,\dots,S_{N-1}$.
But $c_Q$ can occur in $S_j$ ($j<N$) only if $j=N-1$ and $q'(N-2)=Q-1$,
while $c_Q$ can occur in $R_j$ ($j<M$) only if $j=M-1$ and $q(M-1)=Q-1$.
In either case $y_2=\dots=y_Q=x_1$ (since $R_{M-1}$ and $S_{N-1}$ lift
to $\hat L$) and $f_Q$ joins $a$ to $x_1$,
which contradicts Lemma \ref{key}.

\bigskip
Suppose next that $a$ occurs in $V_N$.  Then $p(N)=P$.   The occurrences
of $X$ in $\tilde V_N$ arise as indicated in the proof of Corollary \ref{SNX}.
The relevant subword of $V_N$ has the form:
$$[(x_{p(N-1)+2}b_{p(N-1)+1}^{-1})\dots(x_{P}b_{P-1}^{-1})(\tau(e_P)b_P^{-1})]^{\pm 1},$$
or, if $\gamma$ passes through $a$:
$$[(x_{p(N-1)+2}b_{p(N-1)+1}^{-1})\dots(\tau(e_P)b_{P}^{-1})(x_1a^{-1})\dots(b_{p(1)+1}b_{p(1)}^{-1})]^{\pm 1}.$$
The occurrences of $X$ in $\tilde V_N$
correspond to those $x_j$, $j\ge p(N-1)+2$ in this subword
that are equal to $v$, and also to $\tau(e_P)$ if $\tau(e_P)=v$. 
   If
$a=\tau(e_P)$ then since
$$\tilde S_N\sim a\tilde U_Na^{-1}\tilde V_N\sim XA_1X\dots A_tXW_1X^{-1}B_tX^{-1}\dots B_1X^{-1}W_2$$
we see that the two $a$--letters that occur in the same $W_i$ have the same exponent, and hence both $W_i$ are nontrivial in $G_0$, as required.

\bigskip
If $a=\iota(e_P)$ then $\gamma$ passes through $a$.  Assume for the moment
that $x_1=u$.   Then the unique occurrence of $b_P$ in $\tilde U_N$
lies on the same side of all the $X$--letters as the unique occurrence of
$a$.   Hence the $W_i$ that contains two $a$--letters also contains a
single occurrence of $b_P$.   To prove that this $W_i$ is nontrivial
in $G_0$, it suffices to show that $b_P$ does not occur in any
of the relators $R_1,\dots,R_{M-1}$ or $S_1,\dots,S_{N-1}$ of $G_0$.
But $b_P$ can occur in $S_j$ ($j<N$) only if $j=N-1$ and $p(N-1)=P-1$,
while if $b_P$ occurs in $R_j$ ($j<M$), then $j=M-1$ and $p'(M-2)=P-1$.
In either case $x_1=\dots=x_{P}=u$, contradicting Lemma \ref{key}.

\bigskip
This last argument does not apply if $x_1=v$.   In this case we still
have $x_2=\dots=x_P=u$, and since $a=\iota(e_P)$ it follows from 
Lemma \ref{key} that $\tau(e_P)=v$.

If, say, $W_1=1$ in $G_0$, then $A_t=vb_P^{-1}$ and
$A_tW_1B_t=A_tB_t\ne 1$ in $G_0$, 
since this word contains a single occurrence of $b_P$, which by similar
arguments to the above cannot occur in any of the relators of $G_0$.
Hence no more than one pair of letters $X^{\pm 1}$ in $S_N$ can cancel 
modulo the relators of $G_0$, and so $S_N$, as a word in $G_0*\<X\>$,
definitely involves $X$, as required.

\bigskip
This completes the proof of the Lemma.
\end{proof} 

\begin{cor}\sl
The maps $G_{\pm}\ra G_1$ are injective.
\end{cor}

\begin{proof}
The commutative square

\vspace{15pt}
\hspace{150pt}
\setlength{\unitlength}{0.0075in}%
\begin{picture}(265,139)(245,480)
\thicklines
\put(270,515){\vector( 1, 0){ 90}}
\put(247,600){\vector( 0,-1){ 70}}
\put(382,600){\vector( 0,-1){ 70}}
\put(270,615){\vector( 1, 0){ 90}}
\put(240,510){\makebox(0,0)[lb]{\raisebox{0pt}[0pt][0pt]{$G_-$}}}
\put(240,610){\makebox(0,0)[lb]{\raisebox{0pt}[0pt][0pt]{$G_0$}}}
\put(375,610){\makebox(0,0)[lb]{\raisebox{0pt}[0pt][0pt]{$G_+$}}}
\put(375,510){\makebox(0,0)[lb]{\raisebox{0pt}[0pt][0pt]{$G_1$}}}
\end{picture}

\vspace{-5pt}
\noindent
is a pushout, and the maps $G_0\ra G_{\pm}$ are injective by the lemma.
Hence $G_1$ is the free product of $G_+$ and $G_-$, amalgamated
over $G_0$.
\end{proof}

Let $L_+$ be the 1--complex obtained from $\hat L$ by identifying the
$0$--cells $(0,u)$ and $(0,v)$ to a single $0$--cell $0$.   Then 
$L_+$ is homotopy equivalent to the subcomplex $\hat L\cup X$
of $\tilde L$, and $G_+$ is
a homomorphic image of the free group $\pi_1(\hat L)*\<X\>$, which
is naturally identifiable with $\pi_1(L_+)$.   Let us fix the 0--cell 0
as a base-point for $L_+$, and consider the generating set
$$B_+=\{\theta_e=\tau(e)\lambda(e)^{-1}~;~e\in E(\G)\}$$
for $\pi_1(L_+,0)$.   Note that $B_+$ is not a basis, since the unique
cycle in $T(\G)$ gives rise to a relation $R_0$ among the $\theta_e$.
However, this is the only relation, in the sense that $\pi_1(L_+,0)$
has a one-relator presentation $\<B_+~|~R_0\>$.

Similarly, if $L_-$ is obtained from $\hat L$ by identifying the 
$0$--cells $(1,u)$ and $(1,v)$ to a single $0$--cell $1$, then $G_-$ is 
a homomorphic image of the free group $\pi_1(L_-,1)$, which is generated
by 
$$B_-=\{\phi_e=\lambda(e)^{-1}\iota(e)~;~e\in E(\G)\}$$
modulo a single relator $S_0$ arising from the unique cycle
in $I(\G)$.

\begin{theorem}\sl
The correspondence $\theta_e\leftrightarrow\phi_e$ ($e\in E(\G)$)
induces a group isomorphism $G_+\leftrightarrow G_-$.
\end{theorem}

\begin{proof}
The relation $R_0$ among the generators $B_+$ is precisely the
nullhomotopic path $R_0$ in $L$, which lifts to $L_+$ (indeed to
$\hat L$).   Under the isomorphism $\Psi\co F(B_+)\ra F(B_-)$ induced
by the map $\theta_e\mapsto\phi_E$, this relation $R_0$ is mapped
to $\d_-R_0=R_1$, which is a relation in $G_-$.   Hence we have
an induced homomorphism $\pi_1L_+\ra G_-$.  In order to show that
this in turn induces a homomorphism $G_+\ra G_-$, we must show that
each relation of $G_+$ is mapped to a relation of $G_-$.

Each word $R_i$, $1\le i\le M-1$ is mapped under $\Psi$ to 
$\d_+R_i=R_{i+1}$, which is a relation in $G_-$.   Similarly, for
$1\le j\le N$ we have $\Psi^{-1}(S_{j-1})=\d_-S_{j-1}=S_j$, so
$\Psi(S_j)=S_{j-1}$, which is also a relation in $G_-$.
Hence $\Psi$ induces a group homomorphism $G_+\ra G_-$, as claimed.
Similarly $\Psi^{-1}$ induces a group homomorphism $G_-\ra G_+$,
and these homomorphisms are mutually inverse isomorphisms, by standard
arguments.
\end{proof}

\begin{cor}\sl
$G(\G)$ is isomorphic to an HNN extension of the finitely presented
group $G_1$, with associated subgroups $G_{\pm}$.
\end{cor}

\begin{proof}
This is an easy exercise, given the isomorphism described in the previous
lemma.
\end{proof}

This completes the proof of our main result, Theorem \ref{main}.

\section{Further remarks}\label{more}

In the proof of Theorem \ref{main}, we have relied heavily on
one-relator theory to show that our HNN base $G_1$ is indeed
defined by the relators $R_i$ and $S_i$.  If we look at LOTs
of larger diameter, we no longer have these tools at our disposal.

As long as $I(\G)$ and $T(\G)$ each have only two components (and
hence only one cycle), a great deal of the proof goes through.
Certainly the forward and backward derivatives give rise to
two finite sequences $R_i$ and $S_i$ of relators for $G_1$, but in
order to prove that these relations are sufficient to define $G_1$
we would need to prove a Freiheitssatz for the one-relator
products $(G_0*\<X\>)/S_N$ and $(G_0*\<Y\>)/R_M$.   In our case,
we have used the combinatorics of the diameter 3 situation in a 
nontrivial way to show that $G_0$ is free and that $S_N$ properly
involves $X$ (resp $R_M$ properly involves $Y$) modulo the relations
of $G_0$, from which the Freiheitssatz follows.

\bigskip
It seems reasonable to conjecture in more generality that the 
HNN base $B$ for $G$, generated by $\{xy^{-1},~x,y\in V\}$ will be
finitely presented.   One may construct sets of relations on
this generating set analogous to the $R_i$ and $S_i$ above,
by repeatedly applying the forward derivative construction to
nullhomotopic paths arising from closed paths in $I(\G)$ (analogous
to our $S_0$), and the backward derivative construction to
nullhomotopic paths arising from closed paths in $T(\G)$ (analogous
to our $R_0$).   Provided we restrict attention to simple 
closed paths, only finitely many relations arise in this way, and
one can conjecture that these form a set of defining relators for
$B$.

Before making this conjecture precise, let us first give a geometric
interpretation of these relations.   On the 2--complex $K=K(\G)$ we
define a {\it track} ${\bf T}$ in the sense of Dunwoody \cite{D}
as follows: ${\bf T}$ intersects each 1--cell in a single point,
and each 2--cell in two arcs as in the diagram below.

\vspace{10pt}
\hspace{135pt}
\setlength{\unitlength}{0.0090in}%
\begin{picture}(160,160)(240,440)
\thicklines
\put(240,440){\vector( 1, 0){120}}
\put(240,440){\vector( 0, 1){120}}
\put(400,440){\vector( 0, 1){120}}
\put(240,600){\vector( 1, 0){120}}
\put(240,600){\line( 0,-1){ 80}}
\put(400,440){\line(-1, 0){ 80}}
\put(400,600){\line(-1, 0){ 80}}
\put(400,600){\line( 0,-1){ 80}}
\put(240,520){\line( 1,-1){ 80}}
\put(320,600){\line( 1,-1){ 80}}
\end{picture}

\vspace{10pt}
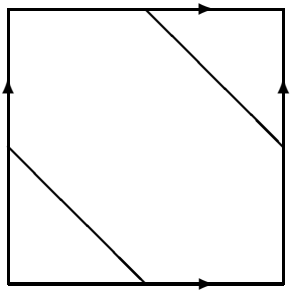
\centerline{\small Figure 2}

\vspace{10pt}
\noindent
The initial graph $I(\G)$ is naturally embedded as a subgraph of the
link of the 0--cell in $K$.   Corresponding to a cycle
$$C=(x_1,\dots,x_n)$$ in $I(\G)$
is a Dehn diagram $D_1$ over ${\cal P}(\G)$ with a single interior vertex 
(whose link maps isomorphically to $C$).   
We also have a nullhomotopic closed path 
$$S_0=(x_1^{-1}x_2)\dots(x_n^{-1}x_1)$$
in $K^{(1)}$.
The boundary label of $D_1$
is $S_1=\d_+S_0$.   Moreover, if we regard $D_1$ as a map from the disc $D^2$ to
$K$, then the track ${\bf T}$ on $K$ induces a track on $D^2$.
This track consists of a single circle in the interior of $D^2$,
together with a collection of arcs, each connecting two adjacent
track points on $\d D^2$.

Now suppose that $S_1$ lifts to $\hat L$.   Then the Dehn diagram $D_1$
can be extended to a diagram $D_2$ with boundary label $S_2=\d_+S_1$,
and so on.   On any Dehn diagram arising in this way, the track
induced by ${\bf T}$ consists of a collection of concentric circles
in the interior of $D^2$, together with a collection of arcs, each connecting two adjacent
track points on $\d D^2$.

Dual to the track ${\bf T}$ is a flow on $K$, indicated on the boundary
of the 2--cells by the arrows in   Figure 2.
The flow induced on $D^2$ by any of the Dehn diagrams obtained as above
has only one singular point in the interior of $D^2$, which is a sink.

\bigskip
We can perform a similar construction for any cycle in $T(\G)$.
The boundary label of the resulting Dehn diagram is obtained by 
repeatedly applying the backward derivative operator to a nullhomotopic
closed path in $K^{(1)}$.
Again, the induced track on $D^2$ consists of a collection of concentric circles
in the interior of $D^2$, together with a collection of arcs, each connecting two adjacent
track points on $\d D^2$.  The induced flow has only one singular point in the interior of $D^2$, which is a source.
 
\bigskip
Let us define a Dehn diagram to be {\it tame} if the induced
track on $D^2$ consists of a collection of concentric circles
in the interior of $D^2$, together with a collection of arcs, each connecting two adjacent
track points on $\d D^2$.   This is equivalent to the induced flow
having only one singular point in the interior of $D^2$, which is
either a sink or a source.   It is not difficult to show that every tame
Dehn diagram arises by the above construction from a cycle in $I(\G)$
or $T(\G)$, and that its boundary label is an alternating word in the
generators $V(\G)$ of $G(\G)$.

\begin{conj}
Let $B$ be the subgroup of $G(\G)$ generated by the alternating words in
$V(\G)$.   Then $B$ has a finite presentation in which the defining relators are the boundary labels of tame Dehn diagrams.
\end{conj}

\Addresses\recd

\end{document}